\documentclass[a4paper,10pt]{article}
\usepackage{amsfonts}
\usepackage{amssymb}
\usepackage{graphicx}
\usepackage{mathtools}
\usepackage{skak}
\usepackage{MnSymbol}

\usepackage{amscd}
\usepackage{MnSymbol}
 
\usepackage{wasysym}

\usepackage{ifsym}

\makeatletter

\newcommand{\Rmnum}[1]{\expandafter\@slowromancap\romannumeral #1@}
\makeatother

\title{  Plateau--Stein Manifolds.}

\author{Misha Gromov}

\begin{document}

\maketitle \tableofcontents

\begin{abstract}  We study/construct (proper and non-proper)  Morse functions  $f$ on complete  Riemannian manifolds $X$
such  that the  hypersurfaces $f(x)=t$ for all $-\infty<t<+\infty$ have   positive mean curvatures at all non-critical points $x\in X$ of $f$.   We show, for instance, that if  an $X$ admits no such (not  necessarily proper)  function, then
   it 
 contains a (possibly, singular)  complete (possibly, compact) minimal 
 hypersurface of finite volume.

\end{abstract}

\section {Introduction.} 

\subsection  {Mean Curvature Convexity.} 

Let $X$ be a smooth Riemannian manifold.   Given a smooth function $f$ on  $X$, define the {\it mean curvatures} $mn.curv_x (f)$
at  {\it non-critical} points $x\in X$, i.e. where $d f(x)\neq 0$, as the mean curvatures of the level hypersurfaces  $Y=Y_r=f^{-1}(r)$, for $r=f(x)$,
$$mn.curv_x (f)=_{def} mn.curv_x (Y_r)$$
where the mean curvatures  of the hypersurfaces  $Y$  are defined by evaluation of their second fundamental forms on the normalized {\it downstream} gradient field
     $-grad(f)/||grad(f)||$.
Call a  function $f$   { \it  strictly mean curvature convex} (sometimes we say  say  "{\sl $(n-1)$-mean curvature convex}"  instead of just "mean curvature convex"  for $n-1=dim (Y) = dim (X)-1$.)   if      $$mn.curv_x (f)  \geq \varepsilon(x)>0   \mbox { for all {\it $f$-non-critical} points   $x \in X$, i.e.   where $df(x)\neq 0$},$$
      for a  positive continuous function  $\varepsilon$ on $X$. 
  
\vspace {1mm}

{\it Remarks.} (a) We are especially concerned with strictly mean curvature convex {\it Morse}  function, i.e.  where the critical points of $f$  are {\it non-degenerate}. Even though our  "convexity" definition {\it formally} makes sense for all smooth  functions $f$,   one has,   in reality, to impose some, possibly weaker than Morse, constrains on the critical points of $f$ ---  we do not want to accept, for example,  constant functions. 

(b) In what follows, most our manifolds  $X$ are non-compact  where    $\varepsilon(x)$ may tend to $0$ for $x\to \infty$. This   happens, for instance,  to the squared distance function  in the Euclidean space  $\mathbb R^n$ from the origin, that is an archetypical example of a   strictly mean curvature convex  Morse  function.
      
(c) If we compose $f :X \to \mathbb R$ with a smooth nowhere locally constant  function
$\psi : \mathbb R\to \mathbb R$, then $\psi\circ f: X\to\mathbb R$ has, at least locally, the same levels as $f$.  The sign of the mean curvatures of the levels is preserved  if  the derivative derivative of $\psi$  is positive,  $\psi'>0$; however, it changes where $\psi'<0$.

\vspace {1mm}

 Conclude by observing that Morse properties of a function are influenced by the sign of the   mean curvature  of the levels
via the following obvious inequality.

\vspace {1mm}

 $[\mathbf {ind{\leq n-2}}]$   {\it the critical points of  strictly mean curvature convex Morse functions have  their Morse indices $\leq n-2$
for $n=dim (X)$. }

\vspace {1mm}

\subsection{  Non-Proper and Proper  Plateau--Stein Manifolds. }  

  \textbf {[n-n-proPS]$_\smile$}.  Call a    possibly non-complete,   Riemannian manifold   $X$ of dimension $n\geq 2$ without boundary  {\it Plateau--Stein}   if it admits a  strictly mean curvature convex  Morse  function $f:X\to \mathbb R$. 
  
  Sometimes, to emphasize that $f$ is {\it not} assumed proper, we call these {\it non-proper}   Plateau--Stein or [n-n-proPS]$_\smile$, where   "$\smile$" stands  for "convex"  with "$\frown$" in the next section for"concave" and where  
  
  \hspace {7mm}  {\it non-proper}  must be always understood as  {\it not necessarily proper.}

  \vspace {1mm}
  
Three other similar conditions on $X$  are as follows.\footnote{Bringing forth these properties was motivated by our conversations with Bruce Kleiner.} \vspace {1mm}

 [PS](1)  Given a compact subset $B\subset X$ and $\varepsilon>0$, there exist a {\it strictly $(n-1)$-volume contracting} continuous map $\Psi_\varepsilon:B\to X$,
 that is
 $$vol_{n-1}(\Psi_\varepsilon(H)) <  vol_{n-1}(H)$$  
  for all smooth hypersurfaces $H\subset X$, and such that  
 $$dist_X(\Psi_\varepsilon(x),x)\leq \varepsilon\mbox { for all $x\in B.$ }$$
  (This condition does not truly need any smooth structure in $X$.) \vspace {1mm} 
 
 [PS](2)  $X$ admits a  $C^1$-smooth {\it strictly $(n-1)$-volume contracting} vector field $V$, $n=dim(X)$, i.e., for every  compact subset $B\subset X$, there exists $\varepsilon>0$, such that  $V$ integrates on $B$ to a  flow   up to  the time $t=\varepsilon$ where the flow maps, say  $V_t:B\to X$,  are  strictly $(n-1)$-volume contracting on  $B$  for  $0<t<\varepsilon$.  \vspace {1mm}

  [PS](3) $X$ admits a {\it strictly $(n-1)$-mean convex function $f$} i.e.   such that the {\it gradient field} of  $-f$ is 
 strictly $(n-1)$-volume contracting.  (See section 3.3 for an alternative definition.)  \vspace {1mm}
 
Clearly,  [PS](3)$\Rightarrow$[PS](2)$\Rightarrow$[PS](1) and also  [PS](3) implies Plateau--Stein,
since  every  strictly $(n-1)$-mean convex function $f$ admits an arbitrarily small perturbation that makes it   strictly $(n-1)$-mean  curvature convex.

  What seems non-obvious   --- I do not see a direct proof of this --- is the following corollary
 to the {\it inverse maximum principle} stated in the next section. 
   \vspace {1mm}

 \textbf{Contraction Corollary for Covering.} {\it Let  $X$ be an infinite covering of a compact manifold. If $X$ 
 is }  [PS](1)  {\it  i.e. if it admits  {\it strictly $(n-1)$-volume contracting} continuous maps $\Psi_\varepsilon:B\to X$, for all compact subsets $B\subset X$ and all $\varepsilon>0$, with $dist_X(\Psi_\varepsilon(x),x)\leq \varepsilon$,
then 
  $X$ is Plateau--Stein.}

   \vspace {1mm}

 {\it Remark/Question.} Probably, if $f$ is strictly mean convex, then one can arrange a smooth  function $a$ on $X$ with a large derivative along the  gradient field $grad(f)$, such that 
 the field $-e^{a(x)}grad(f)$ would be  strictly mean contracting.  
 
 But   it is less clear  what should be {\it exactly} the class of (non-covering)  manifolds  $X$ where
 the existence of a strictly $(n-1)$-volume contracting field implies Plateau--Stein.   
  
  \vspace {2mm}

  \textbf {[proPS]$_\smile$}.  Say that  $X$ is  {\it proper   Plateau--Stein},   if it admits  a {\it proper} {\it positive}   strictly mean convex  Morse function $f:X\to \mathbb R_+$ where proper for (a not necessarily positive) function $f$ on $X$  means that $x\to\infty \Rightarrow |f(x)|  \to\infty $, where 
  "$\to \infty$"  means "eventually leaves every compact subset".

Sometimes we say that {\it the Riemannian metric on $X$ is proper/non-proper  Plateau--Stein}.

 \vspace {1mm}
 
 Proper Plateau--Stein manifolds are reminiscent of {\it complex Stein manifolds}  $X$ that, by definition,   support    proper positive {\it strictly $\mathbb C$-convex}, traditionally called  {\it plurisubharmonic, } functions. An obvious {\it necessary} condition for the existence of such a (not necessarily proper)  function is the absence  of {\it compact complex submanifolds} of positive dimensions in $X$.

A theorem by Grauert   says  that  this condition is also {\it sufficient}, if $X$ can be exhausted by 
compact domain with strictly $\mathbb C$-convex (pseudoconvex) boundaries. 

We shall prove in this paper a Riemannian counterpart to Grauert's theorem with "suitable compact  minimal hypersurfaces" ("hypersurface" always  means a {\it codimension 1 subvariety}, possibly with singularities) instead of   "compact complex submanifolds". 
 
 \vspace {1mm}
 
 The possible topologies of  Plateau--Stein   manifolds are rather transparent. 
 
  \vspace {1mm}
 
   Plateau--Stein manifolds $X$   are {\it non-compact } and,  if  {\it  proper}, they have {\it zero homology}   $H_{n-1}(X;\mathbb Z_2)$  for $n=dim(X)$. In particular, they are   connected at infinity. Moreover,  {\it proper}  Plateau--Stein   manifolds $X$ are diffeomorphic to regular neighborhoods of {\it codimension two} subpolyhedra in $X$.

 \vspace {1mm}
 This follows from $[\mathbf {ind{\leq n-2}}]$.  For example, proper  Plateau--Stein  surfaces are homeomorphic to the $2$-plane $\mathbb R^2$ and proper   Plateau--Stein  $3$-folds are topological handle bodies, 
  while non-proper  Plateau--Stein    allows a complete (warped product) metric on the topological cylinder $X^n=X^{n-1}_0\times \mathbb R$ for all $(n-1)$-manifold $X^{n-1}_0$ as a simple argument shows.

\vspace {1mm}

In fact  one can show  (we leave  this to the reader) the following.  

\vspace {1mm}

{\it Let   $f:
X\to \mathbb R$ be  a proper, not necessarily positive,  Morse function $f:
X\to \mathbb R$, where all  critical points  have  indices $\leq n-2$. Then there exists a complete   Riemannian metric on $X$, (which is proper Plateau--Stein according to our definition } {\it if $f$ is positive) for which this function is strictly mean curvature convex.}

\vspace {1mm}

Unlike  the proper Plateau--Stein the non-proper  Plateau--Stein  condition is not topologically restrictive for open manifolds.

  \vspace {1mm}
 
 {\it every open manifold $X$ admits a  (possibly non-complete) non-proper Plateau--Stein} {\it  Riemannian
 metric.} 
 
 \vspace {1mm}
In fact,   a simple argument shows that
  \vspace {1mm}

 {\it given  a smooth function $f$ without critical points  on a smooth  manifold $X$, there obviously exists a (possibly non-complete)    Riemannian metric on $X$, such that the level hypersurfaces  of $f$ are convex with respect to this function.} 

 \vspace {1mm}

Probably,  every open manifold $X$ of dimension $n\geq 3$  admits a {\it complete} non-proper  Plateau--Stein Riemannian
 metric. 

(Complete  Plateau--Stein   surfaces $X$ 
are
  homeomorphic either  to $\mathbb R^2$
 or to the cylinder $S^1\times \mathbb R$, since  other {\it non-compact complete connected} surfaces necessarily contain (non-contractible)  closed
 geodesics that is incompatible with being
   (complete or not)  Plateau--Stein in dimension $n=2$.)

 \vspace {1mm}
 
 The geometry of  Plateau--Stein  manifolds, unlike their topology, is not  as apparent as  their topology.  
 
  \vspace {1mm}
 {\it  Examples and Questions.} 
 (a)   Let $X = (X, g_0)$ be a  complete simply connected $n$-manifold, $n\geq 2$,   of non-positive  sectional curvature  $\kappa(g_0)\leq 0$. Since  the  spheres $S_{x_0}(R)\subset X$  around a given point $x_0\in X$  
are strictly convex,  such an $X$  is proper Plateau--Stein  and all open subsets $U\subset X$ are  Plateau--Stein.

Even though the inequality $\kappa(g_0)\leq 0$  is unstable  under smooth perturbations of $g_0$ the Plateau--Stein may be stable.  
 
  For example, let the Ricci curvature of $X$  be bounded from below by $-\delta g_0$, $\delta>0$,
 e.g. $X$ is a symmetric space of non-compact type with no  Euclidean factor.
   Then   $mn.curv(S(R))\geq \delta >0$  for all $R>0$.   Since this inequality  is stable 
under uniformly  $C^1$-small perturbations
 $g_\varepsilon$  of the original 
metric $g_0$ on $X$, the function $x\to dist^2_{g_0}(x,x_0)$ remains mean curvature convex with respect to $g_\varepsilon$; hence,    these  $g_\varepsilon$ are proper Plateau--Stein . 

If a  {\it non-flat} symmetric space   $(X, g_0)$ of non-compact type does  have a Euclidean factor, then 
the perturbed metrics  $g_\varepsilon$ are, obviously,  non-proper Plateau--Stein. Probably, they are  
 proper Plateau--Stein.

On the other hand,   the  Euclidean metric $g_0$ on $\mathbb R^n$, $n\geq 2$,  admits arbitrarily $C^\infty$-small     perturbation $g_\varepsilon$ that are {\it not}    Plateau--Stein. 

Indeed, one can arrange $g_\varepsilon$, such that $g_0-g_\varepsilon$ is supported   in an annulus $A$  pinched between two  large spheres, say $S^{n-1}(R)$ and $S^{n-1}(R+1)$  in $\mathbb R^n)$, $R>> 1/\varepsilon $, and such that $(A, g_\varepsilon)$ is isometric to the Riemannian product 
 $S^{n-1}(R)\times [0,1].$ 

It is clear 
such  a $g_\varepsilon$ is
  not  Plateau--Stein, since the mean curvature of a smooth function $f$ on $X=(\mathbb R^n,g_\varepsilon)$ is, obviously, non-positive at the maximum point of $f$ on  $S^{n-1}(R)$.

\vspace {1mm}

(b)  
Let $G_\varepsilon$ be the space of    $\varepsilon$-small    $C^\infty$ perturbations $g_\varepsilon$ of $g_0$ that   
   are  {\it invariant under the action of $\mathbb Z^n$ on}  $\mathbb R^n$, i.e. these  $g_\varepsilon$ correspond to perturbations of the flat metric on the torus $\mathbb R^n/\mathbb Z^n$.  Divide the space $G_\varepsilon$ into
three classes: 
   
   \vspace {1mm}
   
   [1] Proper  Plateau--Stein;
   
   [2] Plateau--Stein  but not  proper Plateau--Stein;
 
   [3]  not even  non-proper Plateau--Stein.
 
 \vspace {1mm}
 
 What is the topological structure of this partition? Are all three subsets  $[1], [2], [3] \subset G_\varepsilon$ dense in  $ G_\varepsilon$  for small $\varepsilon$?
 Is any  of these $[1], [2], [3] \subset G_\varepsilon$  a meager set? (I am uncertain of what happens even for $n=2$ where the answer might be known, albeit in different terms.)

\vspace {1mm}

(c) Let  $p:X_1\to X$  be a {\it Riemannian submersion}  between Riemannian manifolds, i.e.
the differential $dp:T(X_1)\to T(X)$  everywhere has $rank(dp)=n=dim(X)$  and it  is isometric on the  {\it horizontal tangent  (sub)bundle} $T(X_1)\ominus ker(dp)\subset T(X_1)$ (normal the fibers $p^{-1}(x) \subset X_1$). The simplest
instance of this is  the projection of a Riemannian product $X_1=X\times X'$ onto the $X$ factor.

Let the action of the  normal holonomy (by the parallel transport along the horizontal bundle)   on the fibers   be {\it volume preserving}, e.g. $p=X\times X'\to X$. 
Then the $p$-pullback of hypersurfaces from $X$ to $X_1$ preserves their mean curvatures. Therefore,  if a function  $f: X\to \mathbb R$ is Morse strictly mean convex, then a {\it generic} smooth perturbation of  $p\circ f:X_1\to \mathbb R$ is also  Morse as well  strictly mean convex. 

Consequently:

\vspace {1mm}

(c$_1$) {\it If $X$ is  a  Plateau--Stein} {\it    then so is $X_1$}.

\vspace {1mm}

  (c$_2$) {\it  If $X$ is  a proper  Plateau--Stein  and the fibers $p^{-1}(x) \subset X_1$ are closed manifolds for all $x\in X$, then $X_1$ is  also proper Plateau Stein}.

\vspace {1mm}

 (c$_3$) {\it  if  $X$ and  the fibers $p^{-1}(x)$  are proper   Plateau Stein, if  the  action of normal holonomy is  isometric on the fibers, and if the normal holonomy group is compact, then $X_1$ is proper  Plateau Stein.}

\vspace {1mm}

(c$_4$)  It follows from (a) and  (c$_2$) that 

\vspace {1mm}

{\it non-compact   semisimple groups with finite centers are proper Plateau--Stein},

 \vspace {1mm}
 
 while   (c$_1$) implies that 
 
 \vspace {1mm}
 
 {\it unimodular solvable, e.g. nilpotent, groups are (not necessarily proper) Plateau--Stein}.

\vspace {1mm}

 Probably, all   non-compact Lie group $X$ with left invariant metrics, except for  $compact \times\mathbb R$,  are   Plateau--Stein but it is less clear which Lie groups, and Riemannian homogeneous spaces in general,  are {\it proper} Plateau--Stein.

\vspace {1mm}

(d)  The Riemannian cylinders that are product $\mathbb R\times X'$ are  Plateau--Stein  for many (all?) {\it  open} Riemannian $X'$, e.g. for the interiors $X'$  of compact  manifolds with boundaries (this is obvious)  and for complete connected manifolds of with infinite volume. (This is not hard.)

\vspace {1mm}

(e) What  are non-compact,  e.g.   complete, Riemannian manifolds $X'$, such that  the Riemannian products $X\times X'$  are proper  Plateau--Stein for all  proper Plateau--Stein manifolds $X$?

\vspace {1mm}

Conclude with the following questions where
  topology and geometry are intertwined.

\vspace {1mm}
 
Let $V$ be   a  closed  Riemanian manifold.

(I) {\it When  does $V$  admit a Riemannian metric  such that the universal covering $X$ of $V$ with this metric  is proper, or at least non-proper, Plateau--Stein? }

\vspace {1mm}

(II) {\it When is the universal covering of  $V$   {\it bi-Lipschitz  equivalent}
 to a (proper) Plateau--Stein manifold?}\vspace {1mm}

Probably, the answers are invariant under the codimension two surgeries of $V$,  and even  possibly, 
 depend  only on the fundamental group $\Gamma=\pi_1(V)$.



Anyway, the fundamental groups $\Gamma$  of such manifolds  $V$ (where the answers are positive) may be called (I) or (II) "Plateau--Stein"  [proPS]$_\smile$ and [n-n-proPS]$_\smile$.

The best candidates  for "Plateau--Stein" $\Gamma$ are non-amenable groups with one end.
On the other hand,  there may exist some  "tricky" (forget about virtually cyclic)  amenable groups  that are not even   [n-n-proPS]$_\smile$.

 Question (II) makes sense for all complete manifolds $X$, not only coverings of compact ones:
 
{\it When is such an $X$ is   bi-Lipschitz  equivalent to a 
Plateau--Stein manifold?}\vspace {1mm}

("Bi-Lipschitz" seems too restrictive in this context; one needs something half-way from bi-Lipschitz to {\it quasiisometric} in the spirit of the { \it directed Lipschitz metric} \cite {hilbert}.)

\subsection {Inverse Maximum Principle.}
 
 \textbf{[n-n-proPS]$_{\min}$} Say that  a Riemannian manifold $X$  is     [n-n-proPS]$_{\min}$ if it contains {\it  no compact minimal hypersurface. }
 
One can not exclude such a hypersurface being singular. Below is smooth version of this condition with minimal replaced by  "almost concave"

  \textbf{[n-n-proPS]$_\frown$}. Say that  a Riemannian manifold $X$  is     [n-n-proPS]$_\frown$ if it admits a continuos positive function 
 $\varepsilon(x)>0$, such that 
  every  compact smooth domain, i.e. a  relatively compact  open compact subset $U\subset X$ with smooth boundary  in $X$,
 has a point $x\in \partial U$ where $$mn.curv_x(\partial U)\geq \varepsilon(x).$$
 
 Another way to put it is that $X$ contains {\it no}   bounded domain with    $ \varepsilon$-mean-concave boundary. 
\vspace {1mm}

  \textbf{[proPS]$_{\min}$}. Say that  $X$ is  [proPS]$_{\min}$  if it is connected at infinity and if,  for every compact subset $B\subset X$,  there is a (larger) compact  subset $C=C(B) \subset X$, such that 
  
  {\it all compact minimal hypersurfaces 
 $H_i\subset X$   with boundaries contained in  $B\subset X$ are themselves contained in $C$}.

A smooth almost concave version of this condition is as follows.

 \textbf{[proPS]$_\frown$}.  Say that  $X$ is  [proPS]$_\frown$ if it is connected at infinity and  there are  continuos positive function 
 $\varepsilon(x)>0$ and  a proper continuous function $\phi :
 X\to \mathbb R_+$, such that:

given a compact subset $V\subset X$ and a  smooth domain  $U\subset X$ where  
  $$ \sup_{x\in U}\phi(x)>\sup_{x\in V}\phi(x),$$   
there exists a point    $x\in \partial U \setminus V$,  
where $$mn.curv_x(\partial U)\geq \varepsilon(x).$$
  
 \vspace {1mm} 
 
One immediately sees by looking at the {\it maxima points} of  strictly mean curvature convex functions 
$f(x)$ on the boundaries $\partial U$  that proper/non-proper Plateau--Stein manifolds satisfy the corresponding [...]$_{\min}$-conditions. Namely:
  
  \vspace {1mm} 
 
\hspace {29mm}   [n-n-proPS]$_\smile$ $\Rightarrow$    [n-n-proPS]$_\frown$  
 
and

\hspace {36mm}   [proPS]$_\smile$   $\Rightarrow$   [proPS]$_\frown$.
 
  \vspace {1mm} 

 Also it is not hard to see by a simple  approximation argument (see  Step 2 in the next section and  sections 3.4, 5.6, 5.7)
   that     \vspace {1mm}

 \hspace {29mm}   [n-n-proPS]$_\frown$ $\Rightarrow$    [n-n-proPS]$_{\min}$  
 
and

\hspace {36mm}   [proPS]$_\frown$   $\Rightarrow$   [proPS]$_{\min}$.

  \vspace {1mm}

  {\it  IMP for Thick Manifolds.} The main purpose of the present paper is proving inverse implications for  Riemannian
  manifolds $X$  that are {\it thick at infinity} (see section 2.1). Examples of these  include:
  
  $\bullet_{conv}$ \hspace {1mm} complete  manifolds $X$ where the   balls of radii$\leq \varepsilon$ are {\it convex}  for some 
  
  $\varepsilon>0$;
  
  $\bullet_{Ricc}$  \hspace {1mm} complete manifolds $X$ where  $Ricci(X) \geq -const\cdot Riem.metric(X) $ 
  
  and, at  the 
  same time, the  volumes 
  of the unit balls $B_x(1)\subset X$ for all $x\in X$ 
  
  are bounded 
  
  from below by an $\varepsilon>0$;
  
   $\bullet_{Lip}$ \hspace {1mm} complete manifolds   $X$  where the $\varepsilon$-balls  $B=B_x(\varepsilon)\subset X$  for some $\varepsilon>0$ 
   
   admit 
   {\it $\lambda$-bi-Lipschitz} embeddings  $B\to\mathbb R^n$ for $n=dim (X)$ 
   and  for 
   some 
   
   constant $\lambda$
   independent of $B$. \vspace {1mm}
 
 Notice that coverings of compact manifolds are thick at infinity by  either of these conditions.

 \vspace {1mm}

\textbf { Main Theorem:  Inverse Maximum Principle.} {\it  Let $X$ be a  complete     $C^2$-smooth Riemannian $n$-manifold, $n\geq 2$.  If $X$ is thick at infinity,  then}

 \vspace {1mm}
\hspace{-3mm}IMP{[non-proper]}\hspace {14.4mm}     [n-n-proPS]$_{\min}$    $\Rightarrow$   [n-n-proPS]$_\smile$;   \vspace {1mm}

\hspace{-3mm}IMP[proper]\hspace {21mm}  [proPS]$_{\min}$    $\Rightarrow$   [proPS]$_\smile$.

   \vspace {2mm}
 
{\it Convex/Minimal Existence Alternative.} Observe that IMP{[non-proper]} says, in effect,  that either $X$ can be "filled" by {\it strictly  mean convex} hypersurfaces  (that are the levels of a strictly mean curvature convex   Morse function $f(x)$) or, alternatively,  X contains 
 a {\it compact minimal} hypersurface $Y$ and  IMP[proper] encodes a similar alternative.
 
 \vspace{1mm}

 \subsection  { $\phi$-Bubbles,  Plan of the Proof of IMP and Trichotomy Theorem.}

Let  $\mu$ be a Borel measure $\mu$  in  $X$ and define the {\it $\mu$-area} of a {\it domain} $U\subset X$  with boundary $Y=\partial U$   as 
$$area_{-\mu}(U)=_{def} vol_{n-1}(Y)-\mu(U), $$
where "domain" means either a closed subset $U\subset X$ with the interior $int(U)\subset U$ being {\it dense} in $U$ or an open subset that {\it equals the interior} of its closure.

Call $U$  a {\it $\mu$-bubble} 
  if it  {\it locally minimizes}  
 the function $U\mapsto area_{-\mu}(U)$  where  "local"   may be understood at this point 
relative to   the Hausdorff metric in space of pairs $(U,Y=\partial U)$. (We return to this in  section 2.)

 . 
 
 For instance, if $\mu=0$ then $\mu$-bubbles are  domains bounded by stable minimal hypersurfaces in $X$.  

If $\mu$ is given by a measurable density function
 $\phi(x)$,  $x\in X$, i.e.  $ \mu =\phi dx$ for the Riemannian $n$-volume (measure) 
 $dx$, then we speak of    {\it $\phi$-bubbles} and observe that  if $\phi$  is a continuous function, 
 then
 
 { \it the mean curvatures  of the boundaries  
 $Y=\partial U$ of 
 $\phi$-bubbles 
satisfy 

$mn.curv_x(Y)=\phi(y)$ for all regular points $x$ of $Y$.}

\hspace {-6mm} In particular, if $\phi\geq 0$, then $\phi$-bubbles are {\it  mean convex } at all regular points of their boundaries, i.e. $mn.curv_x(Y)\geq 0$ at all regular $x\in Y$ and {\it strictly mean convex} at such points if $\phi > 0$.

In sequel, if $\phi $ is not specified,   these are called  just {\it (strictly) mean convex bubbles }.

(By definition, $Y$ is {\it regular} at a point $x\in Y$ if $Y$ is a {\it $C^2$-smooth} hypersurface in a neighbourhood of this point.)
\vspace {1mm}

{\it Bubbles with Obstacles.}  If $\phi_0 =\infty_B$ equals $+\infty$ on some, say compact subset $B
\subset X$ and zero outside $B$, then the boundaries  of  $\phi$-bubbles $U\supset B$ are, almost by definition,  minimal hypersurfaces  $Y$ in 
the closure of the  complement $X\setminus B$ that solve the Plateau problem 
with the {\it obstacle $B$.}  

We shall often  use  positive continuous functions $\phi=\phi_\varepsilon>0$ that approximate such an  $  \phi_0$, being large, say $1/\varepsilon$,  on $B$ and  $\varepsilon>0$-away from $B$. Then the  corresponding $\phi_\varepsilon$-bubbles $Y=Y_\varepsilon$
lie  close to $\phi_0$-bubbles for small $\varepsilon\to 0 $   that helps to understand their overall geometry, while the continuity of $\phi$ makes the local (quasi)regularity of these $Y_\varepsilon$ similar to that of  
minimal varieties.\vspace {1mm}

We divide the proof of IMP into five steps. \vspace {1mm}

Step 0: {\it Excluding "Narrow Ends".} The representative  case of the theorem  is where $X$ is {\it one-ended}, i.e.  {\it connected at infinity}, and where this end has  {\it infinite area}.  This  means that the boundaries $Y_j=\partial V_j$ of an arbitrary exhaustion of $X$ by bounded domains $V_j\in X$ satisfy
 $$vol_{n-1}(Y_j)\underset{j\to \infty}\to \infty \mbox { for } n=dim(X).$$ \vspace {1mm}
 (In fact, one needs a slightly more general version of this condition as we shall explain in sections 2.2, 2.3.) \vspace {1mm}
  
 Step 1: {\it Mean Convex Exhaustion.} Let  $X$ be one-ended complete with  infinite area at infinity.  Then

  {\it there exist  strictly mean convex compact  bubbles $U_j \subset X$, $j=1,2,...,$   that exhaust $X$,
  $$U_1\subset U_2,\subset ...\subset U_j\subset ...  \mbox   { and } \cup_jU_j=X.$$}     
  
  Notice that  the thickness at infinity condition is not needed beyond this point. On the other hand it is essential for the   existence of  the bubbles $U_j$.
  
 In fact,    such an $\phi_j$-bubble $U_j$ is obtained   for $\phi_j>0$ that is large on the ball $B= B_{x_0}(j) \subset X$  of radius $j$ around a fixed point $x_0\in X$ and small away from $B$.

  The   existence of $U_j$ is proven in sections 2.2, 2.3 by the  standard  minimization argument of the   {\it geometric measure theory}, 
  that works on our (non-comact) $X$   because of 
  the  thickness condition  that is designed exactly for this purpose (see section 2.1)   in order      to prevent  {\it partial escape}   of minimizing sequences of (boundaries of)  bounded domains $U$ in $X$   to infinity.  (Such an  escape can be imagined as  a     "long narrow tentacle" protruding from  the "main body" of $U$).  \vspace {1mm}

 Step 2: {\it Mean Convex Regularization.}  One  can not guaranty at this point that 
the boundary hypersurfaces  $Y_j=\partial U_j$  are smooth for $n=dim(X)\geq 7$. Yet, they do  have positive mean curvatures in a generalized sense.  Moreover,  (this proves the implication [n-n-proPS]$_\frown$$\Rightarrow$
 [n-n-proPS]$_{\min}$) these $Y_j$ can be

\hspace {15mm} {\it approximated by  $C^2$-smooth hypersurfaces $Y'_j\subset U_j$  

\hspace {15mm}  with 
positive mean curvatures. }
 
 In fact, since (the boundaries of) $\phi$-bubbles $U$ are {\it quasiregular}  (as defined in section 3.2)  for all continuous functions $\phi$ by the  {\it Almgren-Allard regularity theorem}, the minus  distance function 
  $d(x)=-dist(x,Y=\partial U)$  can be {\it regularized} almost without   loss of the lower mean curvature bound near the boundary of $U$. Namely, we shall see in section 3.4, and, from a different angle, in section 5.6, that

  \vspace {1mm}

 {\it every strictly mean convex bubble  $U \subset X$ admits a    continuous  
 function
 
 $d': U\to (-\infty,0]$, 
 such that \vspace {1mm}

 $\bullet $  $d'(x)=0 $  for $x\in Y=\partial U,$

 $\bullet $     $d'(x)$ is  smooth strictly negative   in the interior $int(U)\subset U$    
 
 and all critical points  of $d'$ in  $int(U)\subset U$ 
  are non-degenerate,

 $\bullet $ there exists a continuous   function  $\delta$ in $U$ 
 that vanishes on 
 
 the boundary $Y=\partial U$,
 such that $mn.curv_x (d')\geq \phi(x)-\delta(x)$  
 
 for all $d'$-non-critical 
 points 
  $x \in int(U)$.} \vspace {1mm}

 {\it Remarks and Questions.} (a) The above  regularization is non-essential
 at this stage  of the proof; yet,  it will becomes  relevant later on.

 (b) This regularization, along with the simple but "non-elementary"  minimization argument  at Step 1 in the framework of the geometric measure theory   provides  an  exhaustion of $X$ by compact domains $U'_j\subset X$ with smooth strictly mean convex boundaries $Y'_j$. 
 
 Is there   an "elementary" proof of this?

(c) Our  regularization procedures (see sections 3.4, 5.6),   however simple,  needs $C^2$-smootheness of the Riemannian metric in $X$ and  it  does not work for $C^{1}$-manifolds with the sectional curvatures bounded from above and from below.

But the inverse maxima principles, if properly formulated, may hold for $C^1$-smooth manifolds and  for some  singular spaces, e.g. for   {\it Alexandrov spaces} with curvatures bounded from below. 

\vspace {1mm}

Step 3: {\it   Inverse Maximum Principle for Manifolds with Boundary.} Let $V$  be a  smooth compact Riemannian manifold with boundary   and
$\phi$  be a continuous function on $V$ such that the boundary of $V$ is {\it strictly mean $\phi$-convex}, i.e. 
$mn.curv_v(\partial V)>\phi(x)$ for all $v\in \partial V$. Then, assuming $\phi>0$, \vspace {1mm}

$\bullet$ \hspace {1mm}    {\it either the interior of $V$ contains a $\phi$-bubble,

$\bullet \bullet$ \hspace {1mm}   or  $V$ admits a Morse function  
$f: X\to \mathbb R_-=(-\infty,0]$ that vanishes

 on the boundary 
of $X$ and   that is 
strictly  mean $\phi$-convex 

i.e. $mn.curv_v(f)> \phi(v)$  
for all $f$-non-critical points 
$v\in V$. }

\vspace {1mm}
The proof of this is divided into two half-steps. \vspace {1mm}

Half-Step 3A: {\it Shrinking Bubbles.}

  Take a (eventually small), positive $\rho>0$, a (large) $C>1$  and a   monotone decreasing   sequence $\varepsilon_i>0$  $i=0,1, 2,...$, where $\varepsilon _1=\varepsilon_1(V)>0$ is (very) small and where
$\varepsilon _i \to 0$  for  ${i \to \infty}$.

Construct step by step 
a  sequence 
$$U_0=V, \phi_{1} , U_1,\phi_2, U_2,\phi_3,..., U_i,\phi_{i+1}, U_{i+1},... ,$$  
where
$\phi_{i+1}$, $i=0,1,2,...$, is a continuous function on $U_{i}$ and where $U_{i+1}\subset U_{i}$ is a $\phi_{i+1}$-bubble  where the following three conditions must be satisfied.  \vspace {1mm}

{\it $(\ast)_\varepsilon$  \hspace {2.8mm} $\phi_i \geq \phi +\varepsilon_i $ for all $i=1,2,...$; \vspace {1mm}
 
 $(\ast)_{\varepsilon,\rho}$  \hspace {1mm}  $\phi_i(u) = \phi (u)+\varepsilon_i $ for all $u$ in the $\rho$-neughbourhood 
 of $\partial U_{i-1}\subset  U_{i-1}$.   \vspace {1mm}

$({\ast})_{C \rho}$  \hspace {2.8mm}  $\partial U_{i+1}$ is contained in the 
 $C\rho$-neighbourhood of $\partial U_{i} \subset U_i$ 
 
 for all $i=0,1,2,...$.}   \vspace {1mm}

\hspace {-6mm}  This $({\ast})_{C \rho}$  is the only non-trivial requirement on our  sequence, where the existence of a $\phi_{i+1}$-bubble   $U_{i+1} \subset U_i$    with the  boundary  $\partial U_{i+1}$  lying   close to  $\partial U_{i}$   needs a suitably 
chosen  $\phi_{i+1}$  that
must be  large away from the $\rho$-neigbourhood  of 
$\partial U_{i} \subset U_i$ (see sections 2.4,  41,  4.2 ).\vspace {1mm}

Half-Step 3B: {\it Regularization.} On can show (see section 2.)  that if $V$ contain no $\phi$-bubble,
the sets  $U_i$ become empty for large $i$. 

Then, if  $\rho>0$ is sufficiently small,  one can construct the required $f$ by "splicing and regularizing" minus distance functions $u\mapsto -dist(u,\partial U_i)$ on  the subsets $U_i$ (see section 4)\vspace {1mm} 
 
 { \it Remarks.}  (a)    "Shrinking bubbles"  can be seen as  discretization of a 
 some "gradient flow"  for a non-Hilbertian norm in the (tangent space to the) space of subvarieties in $X$, where  the Hilbertian norm leads to the mean curvature flow.

 (b) The condition $\phi>0$ can be dropped with a slightly more  general  notion of "$\phi$-bubble",
 that would allow, for instance, the central geodesic in the (hyperbolic)  M\"obius band for the role of  a $(\phi=0)$-bubble.

Step 4. {\it Limits by Exhaustion.}  We shall use in section 4.2 a simple compactness property  (of sets   of  the boundaries of )  our  mean convex bubbles $U_{i,j} \subset U_j $ construct
mean convex functions $f$ on $X$ as limits of such functions $f_j$ in bounded mean convex  domains $U_j\subset X$  that exhaust $X$.

\vspace {1mm}

Looking closer (see section 4)  one obtains with the above  argument  the following.

 \textbf{   Trichotomy  Theorem}.   Let $X$ be a complete Riemannian $C^2$-smooth manifold (not assumed thick at infinity).   Then  (at least) one of the three  conditions is satisfied.
\vspace {1mm}

(1) {\it    $X$  admits a proper  (positive,  if $X$ is connected at infinity) strictly  mean 

curvature convex Morse function;}\vspace {1mm}

(2) {\it $X$ 
contains a  complete (possibly compact)
minimal hypersurface $H$ of 

finite volume;}\vspace {1mm}

(3)  {\it  $X$  admits a non-proper  strictly  mean curvature 
convex  
 Morse function  
  
  and 
  such 
  that 
  
\hspace {2mm} either there is a non-compact minimal hypersurface $H$ with finite volume  

that 
is closed in $X$ as a subset and that has  compact 
boundary,

\hspace {2mm} or there is a sequence of compact  minimal 
hypersurfaces $H_i\subset X$  
  with no 
 
 uniform   bound on their
diameters,
such that the 
boundaries 
$\partial H_i$ 
are 

contained 
in a 
fixed compact subset
in $ X$.}

\vspace {1mm}

{\it Remarks.}  
(a)  A complete Riemannian manifold $X$ with two ends that
 admits a  proper strictly mean curvature convex  Morse function $f:X\to \mathbb R$,    may contain, however,  arbitrarily large compact  minimal hypersurfaces with boundaries in a given compact subset in $X$. 

For instance, the 2D {\it hyperbolic cusp} $X_0$ (the hyperbolic plane divided  by a parabolic isometry) has this property  and  the Cartesian products $X=X_0\times V$ for compact $V$ furnish example of all dimensions. 

What are other examples of minimal hypersurfaces  protruding toward  "concave ends" in complete manifolds?

Are there such examples with thick ends, e.g. for manifolds with bounded geometries?

(b) The above theorem  (and, in particular, its  special case stated in the abstract to our paper)  shows that the inverse maximum principle does not truly need thickness
at infinity, but   the direct maximum principle, probably, does. Quite likely, there exist  complete
 Plateau--Stein $n$-manifolds  for $n\geq 3$ that contain complete minimal hypersurfaces of finite $(n-1)$-volume.

\subsection {Miscellaneous Remarks, Questions and Corollaries.}

 (A) The most  essential  ingredient of our proof  --- the {\it Almgren-Allard regularity theorem} for "soap bubble" is   trivial  for   $n=dim(X)=2$:  curves with continuous  curvature in surfaces are, obviously, smooth.

 Consequently,  our argument is quite elementary for $n=2$.    In fact,  both IMP  hold  with (almost) no restrictions  on $X$, where  IMP[non-proper] reduces to     the following, most likely known, proposition.

\vspace{1mm}

  \textbf{IMP[dim=2]}.  {\it Let   $X$ be a surface  with a complete $C^2$-smooth} (probably, $C^1$ will do in this case) {\it  Riemannian metric. Then one of the following three properties holds.}  
  
 (1)  {\it  $X$ contains a simple closed geodesic, }

(2) {\it  $X$ supports a smooth function $f:X \to \mathbb R$ which has no critical points and  such that
the sublevels $f^{-1}(-\infty, r]\subset X$ are strictly geodesically  convex  for all $r\in \mathbb R$.} 
(Such an $X$, obviously, is homeomorphic to  the plane $\mathbb R^2$ or to the cylinder $S^1\times \mathbb R$.) 

(3)  {\it $X$ is homeomorphic to the sphere with three point removed.} 

\vspace{1mm}

(Originally I overlooked (3);  it    was pointed to me by Yura Burago that one can not  always ensure a {\it simple} closed geodesic with this topology where non-simpe geodesics are abundant for all metrics. But it should be noticed that the essence  of this   IMP  resides  in the     surfaces that are homeomorphic to $\mathbb R^2$ and to $S^1\times \mathbb R$.)

\vspace{1mm}

\vspace{1mm} 

(B) An instance of a corollary to, say,  IMP[non-proper],  is the validity of the counterparts of the  stability of Plateau--Stein under Riemannian submersions (see  (c) from the previous section) for the corresponding [...]$_{\min}$.

  For example, let  $X$ be a complete  and thick at infinity. \vspace{1mm}

{\it If $X$ is } [n-n-proPS]$_{\frown}$  {\it then so does the Riemannian product  $ X\times X'$, for all closed Riemannian manifolds $X'$.} 
\vspace{1mm} 

This, however, looks almost as  obvious  as  the original   Plateau--Stein  case and, moreover, "thick" seems unnecessary. Indeed, \vspace {1mm}

if a  compact smooth domain
$ U_1\subset X\times X'$ is {\it mean concave}, i.e. its boundary   satisfies 

\hspace {19mm}$mean.curv_x (\partial U_1)\leq 0$ for all $x\in \partial U_1$

  \hspace {-6mm} then  the boundary
of the projection $U \subset X$ of $U_1$ to $X$ is also mean concave  at all 
 {\it regular} points $x \in \partial U$, while singular points have generalized  mean curvatures $=-\infty$. 
This allows a approximation/regularization of $\partial U$ with  $mean.curv (\partial U)\leq \varepsilon$
everywhere for all $\varepsilon>0$.

 \vspace {1mm}
 
On the other hand, the [...]$_{\min}$ counterpart of the above  IMP[non-proper] is not fully trivial.

{\it  If $X\times X'$ contains a closed minimal (possibly singular) hypersurface then so do $X$ and $X'$, provided $X$ and $X'$ are complete  and thick at infinity.}

But the direct  proof of this  by the geometric measure theory is very simple.
 \vspace {1mm}

Notice that minimal hypersurfaces in split  Riemannian manifolds $X_1=X\times X'$ do not always split, e.g.  in flat $3$-tori.  Probably,(?)   there  are {\it non-split}    compact  domains $U_1$ {\it with  minimal boundaries}  in  certain  Riemannian products $X_1=X\times X'$  for  open manifolds $X$ and closed $X'$, where "split" means  $ U_1= U\times X'$. 

 But it seems unclear, for example,  if such  non-split $U_1$ with minimal boundaries  exist in the   products   $X_1=X\times S^1$  of  {\it a  hyperbolic surfaces $X$ of finite areas} by  circles and 
if there are compact domains with minimal boundaries in the products $X_1=X\times X'$ of complete  hyperbolic surfaces $X$ and $X'$ of finite areas. (These $X_1$ are {\it not} uniformly locally contractible but some   IMP may hold.)\vspace {1mm}

(C)  The existence of an exhaustion of  a Riemannian manifold $X$ by compact mean convex domain  is an interesting property in its own right, call it {\it strict mean convexity at infinity}. 

 For instance,  a Galois covering  $X$ of a closed Riemaniinan manifold $\underline X$  is strictly mean convex at infinity unless the Galois  group $\Gamma$ of the covering  is virtually cyclic and  if, moreover, $\Gamma$ is non-amenable,  then $X$  can be  exhausted  by domains with mean curvatures $\geq \varepsilon >0$.
(A representative counterexample for cyclic $\Gamma$ is provided by  manifolds  $\underline X$ that admit  fibrations over the circle, $X\to S^1$, such that the fibers are {\it minimal} hypersurfaces.)

Mean convexity at infinity is visibly  "cheaper" than Plateau--Stein; yet,  there are  non-proper  Plateu-Stein manifolds that are {\it not} mean convex at infinity.  For instance, let  $X$ be a topological cylinder,  i.e. homeomorphic to  $X=X_0\times \mathbb R$ where $X_0$ is a compact manifold and  let 
$g_0$  be a  Riemannian metric on $X_0$.  Let $\phi =\phi(t)$,  $t\in \mathbb R$,  be a positive function and observe that the metric $\phi(t)g_0+dt^2$ on $X$ is {\it concave} at the $t\to -\infty$ end of $X$ rather than mean convex. But such an $X$, obviously is Plateau--Stein since  $f(x_0,t)=\phi(t)$
is a strictly  mean convex function. 

{\it Question}  Under what conditions does  Plateu-Stein  imply  strict mean convexity at infinity? 

(An  easy  instance of such a condition is  {\it thickness at infinity}+  {\it connectedness at infinity}.)

\vspace {1mm}

(D) Smoothness conditions we impose on functions  and on  hypersurfaces  in the  definitions of Plateau--Stein  manifolds and of their   [...]$_{\frown}$-counterparts allows a glib formulation of our  results with no need for concept of "minimal hypersurface".   But insistance on this smoothness 
 look facetious in view of the geometric measure theory  techniques that  underly the essential part of  the argument while   "regularization of bubbles"   that excludes    [...]$_{\min}$  may strike one as a waste of effort. 
 
  In fact, an expected  generalization of the IMP-implications must apply to  non-smooth objects in {\it singular} spaces $X$. On the other hand, the regularization process we employ   delivers --- this is implicit in the arguments in section 4 --- a  simple but non-trivial information on geometry of singular minimal varieties.  (This "information" is by no means new).\vspace {1mm}

(E) Here another obvious consequence of the inverse maximal principle implications 

\hspace {33mm}[...]$_{\min}$$\Rightarrow$[...]$_{\smile}$, \vspace {1mm}

 \hspace {-6mm}  where minimal varieties  and their singularities do not appear. \vspace {1mm}

\textbf{ $(n-1)$-Contraction Corollary.}  
 Let $X$ be a  $C^2$-smooth complete Riemannian manifold that is thick at infinity. \vspace {1mm}

{\it If $X$ admits a strictly $(n-1)$-volume contacting vector field $V$ then 

$X$ is Plateau--Stein. \vspace {1mm}

If, moreover, $X$ is connected at infinity and if there are vector fields $V_i$, 

such that the supports of $V_i$ exhaust 
$X$ and such that  $V_i$ are strictly 

$(n-1)$-volume contacting in the 
complement 
of their supports,
then 

$X$ is proper Plateau--Stein.} \vspace {1mm}

(Recall, that  this theorem  was stated for coverings of compact manifolds in section 1.2   and that  a vector field $V$ is {\it strictly $(n-1)$-volume contacting} if  the $V$-derivatives of the volumes of all smooth hypersurfaces  in $X$  are negative.) \vspace {1mm}

(F) Let  $X$, not necessarily thick at infinity,  contain  no {\it compact} minimal hypersurface. Then does it  admit a strictly $(n-1)$-volume contracting vector field? 

(This question is motivated by such a result for $1$-volume (i.e.  length) contracting fields that was pointed out to me by Bruce Kleiner.  Possibly, there is something like that for all  {\it $k$-volume contracting} vector fields.)

\section { Construction of $\phi$-Bubbles and of Minimal Hypersurfaces.}

We shall describe  in this section  a few  standard $\phi$-area (including  (n-1)-volume) minimization constructions   that  deliver  minimal hypersurfaces, such as 
 $\phi$-bubbles, under the thickness condition.

\subsection {Thickness at Infinity.}

 An $n$-dimensional  Riemannian manifold   $X$ is called 
{\it thick at infinity} if it contain {\it no non-compact  minimal}  hypersurface with {\it compact} boundary and with   {\it finite} $(n-1)$-volume. \vspace {1mm}

 Such a hypersurface $Y\subset X$  must, by definition, be closed in $X$  as a subset and be  {\it $\varepsilon$-locally   $vol_{n-1}$-minimizing in $X$  at infinity.}
 This means that  \vspace {1mm}

    there exists an $\varepsilon >0$ ($\varepsilon =1$ is good for us) and  a compact subset  $A=A(Y)\subset X$ (that contains the boundary of $Y$) such that the intersection $Y\cap B$ with 
 every $\varepsilon$-ball $B=B_x(\varepsilon)\subset X$ for $x\in X\setminus A$ is $vol_{n-1}$-minimizing in $X$ in the class of hypersurfaces (integral currents) with the boundary equal $Y\cap\partial B$.  \vspace {1mm}

{\it Examples.} The  Paul L\'evy (Buyalo-Heintze-Karcher) tube volume bound  shows that the  condition $\bullet_{Ricc}$ from section 1.3 implies this thickness, while the conditions  $\bullet{conv}$ and $\bullet _{Lip}$ are taken care by the following corollary to the implication 

\hspace {20mm}[$cone$ $inequality]\Rightarrow$ [$filling$ $inequality$] 

\hspace {-6mm} and the  lower bound on volumes of minimal varieties by the  filling constant     \cite{wenger}.\vspace{1mm}

$\bullet_{fill}$  {\it If every  closed integral $k$-chain  $S$ in $X$, $k=1,2,...$, of diameter $\leq 1$ bounds a $(k+1)$-chain $T$ such that 
$$vol_{k+1}(T)\leq const \cdot diam(S)\cdot vol_k(S)$$
for some $const=const(X)$, then $X$ is thick at infinity.}\vspace{1mm}

 {\it Remark.} This thickness  concept obviously generalizes to all dimensions $2\leq m\leq n-1$  with minimal subvarieties $Y\subset X$ of dimension $m$, where the above criterion remans valid, while   $\bullet_{Ricc}$ should be replaced by $\bullet_{sect}$ with the lower sectional curvature bound instead
 of  such  a bound on  Ricci.

\subsection {Convex and Concave Ends.}  Let $X$ be a  Riemannain manifold possibly with compact boundary $\partial X$
and a single  end, such that $X$ is {\it thick at infinity}. Then one of the following three possibilities is realized. \vspace {1mm}

 [1$_\smile$]  {\it $X$ can be exhausted by compact  strictly mean convex bubbles};
 
 [2$_\frown$]   {\it $X$ can be exhausted by compact  strictly mean concave bubbles};
 
 [3$_\equiv$] {\it There exists a continuous positive  proper function  $h: X\to \mathbb R_+,$
 
 such that 
 the  levels $Y_t= h^{-1}(t) \subset X$ are minimal hypersurfaces, 
 
 that are the  boundaries of $0$-bubbles,
 for all $t\geq t_0=t_0(X)\geq 0$.} \vspace {1mm}

{\it Proof.} Start by observing that 1  and 2 are {\it not} mutually exclusive and if there are these two kinds of exhaustions then there also (obviously) exists an exhaustion by $0$-bubbles. 
But 3  is incompatible with 1    and with 2  by the maximum principle.  

Let  3 do not hold and, moreover, assume that  $X$ {\it can not be exhausted} by compact  $0$-bubbles.
Then minimization of  $vol_{n-1}(Y)$ for $Y=\partial U$ where $U\subset X$ is a compact sufficiently  large domain, either move $Y$ to infinity or brings to a compact region $X_0\subset X$.

In the former case, let $-\phi$ be a negative  function on $X$, where $\phi$  is very large at infinity and very small
in the vicinity of some $X_0$ that contains the boundary of $X$.  Then minimization of $\phi$-area
brings $U$ to a strictly mean concave   $-\phi$-bubble in $X$. 
 
Similarly, if volume   minimization brings $Y$ to a fixed compact $X_1\subset $  we use a positive $\phi$ that is very large 
 on some compact  $X'_1\supset X_1$ and very small at infinity. Thus we obtain a   strictly mean convex bubble.
 
And, keep   obviously modifying $\phi$ by moving the change of its size to infinity,  we exhaust $X$ either by strictly mean convex or by strictly mean concave bubbles.
 
Now, let $X$ be exhausted by  by compact  $0$-bubbles. Then $X$ contain infinitely many "empty bands", say $W$ between  the boundaries of these bubbles  say $U_1$ and $U_2\supset U_1$, i.e. $W=U_2\setminus U_1$, where such a $W$ is  bounded by  the minimal hypersurfaces $Y_1=\partial U_1$ and $Y_2=\partial U_2$ with no $0$-bubble between the two. Then the obvious adjustment of the above argument delivers both a strictly mean convex and a strictly mean concave bubble pinched between $U_1$ and $U_2$.

 Besides,  {\it Almgren's min-max argument}  delivers a non-stable minimal hypersurface in $W$   that separates $Y_1$ from $Y_2$.

\subsection {Minimal Separation of Ends.}  Here $X$   has several ends and  no boundary,
where the set of ends is given its natural topology. 

Notice that every {\it isolated} end $E$ can be represented/isolated  by an equidimensional  submanifold (domain)  $X_E\subset X$ with  compact boundary and a single end;  we say "exhaustion of $E$" instead of "exhaustion of $X_E$ for some $X_E\subset X$".

If the set of ends of  $X$ 
contains at least two limit points, (e.g. if it has no isolated ends),  then $X$ contains a compact two-sided  
smooth hypersurface  $H\subset X$ such that there are infinitely many ends of $X$ on either side of $H$. 
Then, clearly,\vspace {1mm}

 {\it if $X$ is complete and thick at infinity, then  the homology  class of  $H$ contains a minimizing hypersurface.} 

\vspace {1mm}

Now let $X$ has { \it at least two isolated}  ends, say $E_1$ and $E_2$. 

\vspace {1mm}

{\it If none of these ends admits a strictly mean  convex exhaustion, then $H$ contains  a compact (non-stable)  minimal variety $M\subset X$ where this $M$ may be a varifold.} \vspace {1mm}

{\it Proof.}  Let $h:x\mapsto t=h(x) \in (-\infty,+\infty)$ be a smooth proper Morse function $X \to \mathbb R$ such that $h(x)\to  -\infty$ for $x\to E_1$   and $h(x)\to  +\infty$ for $x\to E_2$.
 
 The $vol_{n-1}$-minimimzation process starting from the levels
$h^{-1}(t)$  moves some connected  component of $h^{-1}(t)$ for small negative  $t$ (approaching $-\infty$) to $E_1$,
while some component for positive  large $t$ goes to $E_2$.  Since the manifold $X$ is thick at infinity,  Almgen's min-max theorem applies and the proof follows. \vspace{1mm}

By combining the above with  [1$_\smile$]/[2$_\frown$]/[3$_\equiv$] in the previous section, we conclude to the following.

 Let $X$ be a complete Riemannian manifold that is thick at infinity. Then  
\vspace{1mm}

 $[\smile/-]$  {\it Either  an isolated end of $X$  admits a strictly mean convex exhaustion,

or  $X$ contains a compact minimal subvariety.}

\subsection {Shrinking Mean Convex Ends.}

 Let $X$ be a complete connected Riemannian manifold of dimension $n$ with non-empty compact boundary $Y_\smile=\partial X$ with $mn.curv(Y_\smile)>\varepsilon_0>0$, let $0<\varepsilon _i< \varepsilon_0 $, $i=1,2,... $ be a sequence of positive numbers that converges to 0 and let $\rho(x)>0$  be a continuous function on $X$.  Then  either \vspace {1mm}

\hspace {-3mm} $\star_A$ \hspace {1mm}  {\it  $X$ contains a  minimal hypersurface $H\subset X$ of finite $(n-1)$-volume 

that is closed in $X$ as a subset and that does not meet $Y_\smile$}, \vspace {0mm}

or

\hspace {-3mm} $\star_B$ \hspace {1mm} {\it $X$ can be  exhausted by  an increasing   sequence of compact strictly concave bubbles $U_i$  in $X$ that contain $Y_\smile$,
$$Y_\smile\subset  U_1\subset U_2 \subset...\subset U_i\subset... \subset X,  \mbox {  }  \cup_i U_i=X,$$
such that 
   
 $\bullet$   $U_{i}$ is contained in the $\rho_{i-1}$-neighborhood of $U_{i-1}$
 $$\mbox { for }  \rho_{i-1}=\inf_{x\in U_{i-1}}\rho(x)\mbox {  and   all $i=1,2,...$,  where  $U_0=_{def}Y_\smile$};$$,

 $\bullet$  the mean curvatures of the topological boundaries $Y_i=\partial U_i \subset X$ satisfy
    
  \hspace {2mm}  $mn.curv_x(Y_i)=-\varepsilon_i $  at all those regular     $x\in Y_i$ 
  where $dist (x, U_{i-1})\leq \rho/2$.} \vspace {1mm}
 
 {\it Remarks}.  
 If  $X$ is thick at infinity, then, as we know, the above minimal hypersurface  $H$, if it exists at all, must be compact. 
 
  If $X$ is compact and no minimal $H$ exists, then the sequence $U_i$ stabilizes and the boundaries$Y_i$ of the bubbles $U_i$ become  empty for large $i$. On the other hand if $H$ does exist and $\rho\leq \rho_0=\rho_0(H)>0$,  then {\it no $U_i$  intersects} $H$. (If $\rho $ is large  then $U_i$ {\it may} exhaust $X$ even in a presence of $H$.) \vspace {1mm}

 {\it Proof of  the $\star_A/\star_B$-Alterantive.} Proceed as at the half-step 3A from 1.4. Namely,  granted  $U_{i-1}$  for some $i$,  let  $\phi_i$  be a  positive continuous  functions   $\phi_i$ on $X$,  such that
 
  $\phi_i =\varepsilon _i$ in the $\rho_i/2$-neighbourhood of $U_{i-1}$,
 
$ \phi_i(x)$ is very large for $dist(x, U_{i-1})\geq \frac {2}{3}\rho_i.$
 
 Then, clearly, there exists a  compact $-\phi_i$-bubble $U_i \supset U_{i-1} $ that satisfies all of the above properties. 
 
 Since $\varepsilon_i\to 0$ the boundary of the union   $\cup_iU_i \subset X$ must be a {\it minimal}
 hypersurface  $H$ in $X$   with $vol_{n-1}(H)< vol_{n-1}(Y_0)$; if no such hypersurface exists,    then this boundary must be empty and  $\cup_iU_i= X$. QED. \vspace {1mm}

 If $\star_B$ holds for all $\rho(x)>0$ then, by the maximum principle, $X$ contains no minimal hypersurface $H$ of {\it any volume}.  This leads to the following\vspace {1mm}

\textbf {Almgren's min/max  Theorem for  Non-Compact Manifolds.} {\it If a complete Riemannian  manifold $X$ with non-empty compact 
strictly mean convex boundary contains a complete minimal hypersurface  then it also contains a  complete minimal hypersurface of finite volume.} 

("Complete" means  being closed in $X$ as a subset,  not intersecting the boundary of $X$ and having no boundary of its own.)

\subsection {Shrinking to Concave Boundary.} The above admits a relative version where 
$X$ has a  concave component  $Y_\frown$   in its boundary (or several such components) that serves as an obstacle for shrinking bubbles and where
the minimal hypersurface $H$ that (if it exists at all) obstructs shrinking of bubbles may have non-empty  boundary that is contained  in $Y_\frown$.

Namely, let   again $X$ be a complete connected Riemannian manifold of dimension $n$ with  compact boundary that is now decomposed into two disjoint parts 
$\partial X=Y_\smile\cup Y_\frown$ (these $Y_\smile$ and  $Y_\frown$ are unions of connected components  of $\partial X$) where  $Y_\frown$ (possibly, empty as in the previous section)  is strictly mean concave and where $Y_\smile$ is non-empty and has 
$mn.curv(Y_\smile)>\varepsilon_0>0$. 

Let $\rho(x)>0$  be a continuous function on $X$.  Then  either \vspace {1mm}

\hspace {-3mm} $\frown\star_A$ \hspace {1mm}  {\it  $X$ contains a  minimal hypersurface $H\subset X$ of finite $(n-1)$-volume

that is closed in $X$ as a subset, that does not meet $Y_\smile$, 
and that 

 may have 
boundary contained in  $Y_\frown$, \vspace {0mm}}

or

\hspace {-3mm} $\frown\star_B$ \hspace {1mm} {\it the complement $X\setminus Y_\frown$ can be  exhausted by  an increasing   sequence of compact strictly concave bubbles $U_i$  in $X$ that contain $Y_\smile$,
$$Y_\smile\subset  U_1\subset U_2 \subset...\subset U_i\subset... \subset X,  \mbox {  }  \cup_i U_i=X\setminus Y_{\frown},$$
such that 
   
 $\bullet$   $U_{i}$ is contained in the $\rho_{i-1}$-neighborhood of $U_{i-1}$
 $$\mbox { for }  \rho_{i-1}=\inf_{x\in U_{i-1}}\rho(x)\mbox {  and   all $i=1,2,...$,  where  $U_0=_{def}Y_\smile$}.$$}
 
  {\it Remarks}.   We could also  impose here additional constrains on the mean curvatures of 
 the  boundaries $Y_i=\partial U_i \subset X$ as we did in the absence of $Y_\frown$ but we postpone doing this until section 4 where this becomes relevant.

 An essential case of the above is where $X$ is compact with two boundary components, $Y_\smile$ and $Y_\frown$ and these shrinking bubbles are  used (see section ....) for construction of a strictly mean curvature convex Morse function  $f$ on $X$ that equals $1$ on $Y_\smile$  and $0$ on
 $Y_\frown$, where

 {\it such an $f$ exists if and only if there is no minimal hypersurface 
 
 $H$ in $X$
 with $\partial H\subset Y_\frown$.} 
 
\hspace {-6mm}  Here (and  everywhere in this kind of context) "only if" follows by the maximum principle, while "if" is what we call   the inverse maximum principle.\vspace {1mm}
  
 {\it Proof of $\frown\star_A/\frown\star_B$. }  Proceed as earlier and keep  pushing boundaries of  bubbles closer and closer to $Y_\frown$. Then, say in the compact case, we arrive at  a  {\it maximal}  compact bubble $U_{max}\subset X$ the boundary of which can  not be moved  closer to $Y_\frown$ anymore.  Then  either  $U_{max}=X$  or the topological boundary of  $U_{max}$ is non-empty.  Then this boundary, call it $H$ makes our  {\it minimal} hypersurface  in the interior of $X$: if not "minimal" it could  be moved closer to $Y_\frown $.
 
 Notice that  this $H$ is {\it tangent} (rather than transversal) to $Y_\frown$ where the two hypersurfaces meet.

 \section {Distance Functions, Equividistant Hypersurfaces and  $k$-Mean Convexity.}

  We fix in this section our terminology/notation and  state   a few   standard facts on distance functions  in Riemannian manifolds $X$.

   \subsection{Signed Distance Function and Equidistant Hypersurfaces.}

\vspace {1mm}

 {\it Interior Domains  $U^{<}_{-\rho}$, $U_{-\rho}=U^{\leq}_{-\rho}$ and Equidistant Hypersurfaces  $Y_{-\rho}$.} Let  $U$ be a (closed or open)  domain or an open subset (possibly dense)    in 
   a Riemannian manifold $X$ (or in any metric space for this matter) and    denote by  $ x\mapsto d(x)= d_U(x) = dist_\pm(x, \partial U)$ the {\it signed distance function} to the topological boundary  $Y= \partial U $, i.e.    
   
   \vspace {1mm}
   
  \hspace {10mm} $d_U(x)$ equals the distance from $x$ to  $Y$ outside $ U$,

   \hspace {10mm}    $d_U(x)$  equals {\it minus} the distance from $x$ to  $Y$ in $U$. 
   
   In writing,
   $$d_U(x)=dist (x,Y)=_{def}\inf dist_{y\in Y} (x,y) \mbox { for all } x\in X\setminus U,$$
   $$ d_U(x)=-dist (x,Y)=-dist (x,X\setminus U)\mbox  { for all } x\in U,$$
where
 $$dist (x,Y=\partial U)=dist(x,U), \mbox { for $x\in X\setminus U$},$$
 since the Riemannian distance  is a {\it length metric} being defined via the lengths of curves between pairs of points;   
   \vspace {1mm}

Let $\rho\geq 0$ and denote by
     $U_{-\rho}= U_{-\rho}^\leq\subset U$  and   $U_{-\rho}^<  \subset U_{-\rho}  $, $    \rho\geq 0$, the  closed/open  $(-\rho)$-sublevels of $d_U$, that are
$$U^\leq_{-\rho}=d_U^{-1}(-\infty, -\rho] \mbox { and }  U_{-\rho}^<=d_U^{-1}(-\infty, -\rho), $$  
  where, clearly,   $U_0=U_{-0}=U^{\leq}_{0}$ equals the topological  closure  of $ U$ and  $U^{\leq}_{-\rho}$ are closed subsets in  $U$  for $\rho>0$ with  $U_{-\rho}^<$ being  equal the interior of  $    U_{-\rho}^\leq   $ for $\rho>0$.

Let 
 $$Y_{-\rho}=d_U^{-1}(-\rho)\subset   X, \rho\geq 0,$$
be the {\it (interior) $\rho$-equidistant hypersurface to Y},  that is the  subset of points $u$ in  $U$, where     $dist_X(u, Y=\partial U) =\rho$ and that equals  the topological boundary
$\partial U^{}_{-\rho}$ since Riemaniinan manifolds $X$ are {\it  length metric spaces}.

Similarly define  $U_{+\rho}= U_{+\rho}^\leq$  and   $U_{+\rho}^< \subset    U_{+\rho}^\leq   $ or $\rho\geq 0$ as 
$$U^\leq_{+\rho}=d_U^{-1}(-\infty, \rho] \mbox { and }  U_{+\rho}^<=d_U^{-1}(-\infty, \rho). $$  
Thus, $U_{+\rho}= U^\leq_{+\rho}$ equals the {\it closed $\rho $-neighbourhood of} $U$  in $X$ and   $U_{+\rho}^<$ is the {\it open}  $\rho $-neighbourhood.

\vspace {1mm}

{\it On Hausdorff (Dis)Continuity.} Clearly, the boundaries of the open sublevels of $d_U$ satisfy
$$\partial  U^{<}_{-\rho} \subset Y_{-\rho}=\partial  U^{\leq}_{-\rho},$$
 where the {\it  local minima }  of the $d_U$ on $U$ make the difference set   $Y_{-\rho}\setminus \partial  U^{<}_{-\rho} $.

The set valued function  $\rho\mapsto  U^{}_{-\rho}\subset X$, $\rho \in \mathbb R_+$, is  continuous    for the Hausdorff metric in the space of subsets in $X$ at those $\rho$ where   $\partial  U^{<}_{-\rho} = Y_{-\rho}$, or, equivalently, where the closure of the interior of   $U_{-\rho}$ equals $U_{-\rho}$.
 Since $\rho\mapsto  U^{}_{-\rho}$ is a  {\it monotone decreasing} function in $\rho$ for the inclusion order on subsets,  it has {\it at most countably many}
discontinuity points $\rho$.

 Also observe that the function   $\rho\mapsto  Y_{-\rho}=d_U^{-1}(-\rho)=\partial  U^{}_{-\rho}$ is Hausdorff continuous at the Hausdorff 
 continuity points  of the function   $\rho\mapsto  U^{}_{-\rho}$ and 
the word  "hypersurface" is  justifiably applicable to $Y_{-\rho}$ at these  continuity points $\rho$. \vspace {1mm}

{\it Exercise.} Let $Z\subset $X be a compact subset that is contained in a smooth hypersurface in $X$.
Then, for all sufficiently small $\rho>0$,  there exists an open subset   $U\supset Z$  in $X$ with smooth boundary $Y$, such that $Z=Y_{- \rho}=U_{- \rho}$, i.e.  $Z$ serves in $U$ as  the  set of the minima of the (minus distance to $Y$) function $d_U:U\to (-\infty, 0)$.

\vspace {1mm}

  {\it  Example: $U^{}_{-\rho}$ as the Intersection of Translates of $U$.}  If $X=\mathbb R^n$  then, obviously, 
$U^{}_{-\rho}$ equals the intersection
of the  parallel $r$-translates $U+r\subset  \mathbb R^n$ for all $r\in  \mathbb R^n$ with $||r||\leq \rho$ and $Y_{-\rho}$  equals the topological boundary of this intersection.
  $$U^{}_{-\rho}=\bigcap_{||r||\leq \rho} U+r \mbox { and } Y_{-\rho}= \partial (\bigcap_{||r||\leq \rho} U+r). $$
Thus,

{\it the transformation $U\mapsto U_{-\rho}$ 
 preserves all  classes of  Euclidean
 domains 

 (e.g. the class of    convex domains) that are closed under intersections. }

\hspace {-6mm } It is also clear that 
 $$Y_{-\rho}\subset \partial (\bigcap_{||r||= \rho} U+r)$$
and
if the boundary of $U$ is {\it connected}, then 
 $$U^{}_{-\rho}=\bigcap_{||r||= \rho} U+r \mbox { and } Y_{-\rho}= \partial (\bigcap_{||r||= \rho} U+r).$$

More generally, let   $iso_{\leq\rho}$ denote the set of isometries $r: X\to X$ such that $dist(x, r(x))\leq \rho$  for all $x\in X$. Then, obviously,
$$U^{}_{-\rho}\subset \bigcap_{r \in  iso_{\leq \rho}} r(U).$$
Furthermore, if 
$X$ is a {\it compact  two-point homogeneous space}, i.e. the isometry group of $X$ is transitive on the unit tangent bundle of $X$, then, as in the Euclidean case, 
$$U^{}_{-\rho}= \bigcap_{r \in  iso_{\leq \rho}} r(U)\mbox { and } Y_{-\rho}\subset \partial (\bigcap_{r \in  iso_{= \rho}}  r (U))\mbox {  for }   iso_{= \rho}=\partial  ( iso_{\leq \rho}).$$
 
{\vspace {1mm}

 \vspace {1mm}   
\subsection { Accessibility and Quasi-regularity.}

A point $x$ in the boundary   $Y=\partial U$, is called  {\it $\rho$-accessible} (from  $U$),
  $$x\in(U_{\rho})_{+\rho}\subset U\cup Y.$$
In other word, $x$ is contained in some Riemannian $\rho$-ball in $X$ that is contained in the closure of  $U$. (The referee pointed out to me that this  is   usually called "with reach $\rho$", with a possible origin of the concept due to
 Federer.)

Say that  an open subset  $U$  in $X$ is {\it $C^2$-quasiregular} (at its boundary)    if, loosely speaking, the singular locus    $ sing_Y\subset Y=\partial U$  is non-accessible from $U$.
   More precisely,  the following two conditions must be satisfied.\vspace {1mm}

(1)   The subsets 
$$Acc_{>\rho_0}(Y)=A_{>\rho_0}(Y)=\bigcup_{\rho>\rho_0}A_\rho(Y) \subset  Y\subset X$$  
  are {\it open} in $Y$ for all $\rho_0\geq 0$.
  
  Notice that this condition  implies that, besides being open, the subsets 
$A_{>0}(Y_{-\rho}) \subset Y_{-\rho} $, $\rho>0$,
  are $C^{1,1}$ in $C^{1,1}$-smooth Riemannian submanifolds in $X$.

(2)  The subset  
$A_{>0}(Y) \subset X$
 is 
 a  {\it $C^2$-smooth} hypersurface in  $X$, that is a  $C^2$-smooth $(n-1)$-submanifold without boundary that, topologically, is a locally closed subset in $X$.\vspace {1mm}

If $U\subset X$ is a {\it closed}  domain then  its  quasiregularity  means that of the {\it interior} $int(U)\subset X$. 

On the other hand, "quasiregularity of a hypersurface" $H\subset X$  is understood as 
 quasiregularity if its complement 
$X\setminus H\subset X$.

\vspace {1mm}

{ \sc Almgren-Allard Quasiregularirty Theorem}. { \it Let $X$ be a    $C^2$-smooth Riemannian manifold.
Then $\phi$-Bubbles $U\subset X$ are $C^2$-quasiregular for all continuous functions $\phi(x)$. Also, all kinds of  minimal hypersurfaces $H\subset X$ are quasiregular.} (See  \cite{white} for a simple prove of this.)

\vspace {1mm}

 Here  "minimal hypersurface" is understood as  a {\it minimal varifold} that does not, necessarily, bound any domain in $X$.

\vspace {1mm}

The following two instances of quaisiregularity are, unlike Allmgren-Allard theorem, are fully obvious.   \vspace {1mm}

(A)  {\it Locally finite  intersections of $C^2$-smooth domains with transversally 
 
 intersecting boundaries are  $C^2$-quasiregular.}

And

(B)  {\it If  $U$ a  $C^2$-quasiregular domain in a $C^2$-smooth Riemannian manifold, 

e.g. 
the boundary $Y$ is $C^2$-smooth to start with, 
then
the sub-domains $U^<_{-\rho}$ 

are also $C^2$-quasiregular.}

 \vspace {1mm}

\subsection{ Smooth and non-Smooth $k$-Mean Convex Functions and Hypersurfaces.}

Let  $X$ be a $C^2$-smooth Riemannian manifold and let $Gr_k(X)$ be the Grassmann space of the tangent 
  $k$-planes   $\tau$ in $X$. Define the {\it $k$-Laplacian $\Delta_k$} from $C^2$-functions $f$ on $X$ to functions on $Gr_k(X)$ by taking the traces of the Hessian of $f$ on all  $\tau\in Gr_k(X)$ 
  $$ \Delta_k(f)(\tau)=trace_\tau Hess(f).$$
For example $\Delta_n=\Delta$ is the ordinary Laplacian for $n=dim(X)$ and $\Delta_1(f)(\tau)$ equals the second derivative of $f$ on the geodesic in $X$ tangent to the  tangent line $\tau$.

Say that $f$  is { \it $k$-mean $\varphi$-convex} for a given continuous function $\varphi=\varphi(\tau)$ on  $Gr_k(X)$ if  $\Delta_k(f)(\tau)\geq $, where, as usual strictly corresponds to $>\varphi $  and plain "convex " stands for $0$-convex.

Observe that the  $k$-mean convexity says, in effect that the  gradient of $-f$ is strictly $k$-volume contracting.   Also notice that
$$ \mbox {$k$-$mean $ $convex$ $\Rightarrow $$l$-$mean $ $convex$ for $l\geq k$},$$
 $$\mbox {$1$-$mean $ $convex$ =$convex$},  \mbox { \hspace {1mm}n-$mean $ $convex$=$subharmonic$}$$    
and that 
$$ (n-1)\mbox{-{\it mean convex}} \Rightarrow \mbox {  {\it mean curvature convex}, }$$
while the converse implication is not, in general true.

On the other hand,  a $C^2$-smooth strictly mean convex  co-oriented hypersurface $Y\subset X$ (e.g.  if $Y=\partial U$)  can be realized as

{\it the zero set of 
a $C^2$-smooth strictly $(n-1)$-convex function $f(x)$ defined  

in some neighbourhood of $Y$.}\vspace {1mm}

{\it Bending $d(x)$ to  a Mean Convex $f(x)$.} (Compare \cite {lawson}, \cite {sha}.)  A strictly mean convex  function $f$  can be obtained, for example, by "bending" the  signed distance function $d(x) =\pm dist(x,Y)$ (that is $d(x)=d_U(x)$ for $Y=\partial U$), i.e. where  bending is achieved by means of     a smooth strictly monotone increasing function $\beta(d)$, $ -\infty <d<+\infty $,  that vanishes at $d=0$, that has the first derivative $d'(0)=1$ and  positive second derivative $d''(x)>0$. 

If $Y$ is compact, then the function
 $$f(x)=\beta\circ d(x)= \beta(d(x))$$  
is strictly $(n-1)$-mean convex  in some neighborhood of $Y$, provided the second derivative $d''(0)$ is sufficiently large (compared to the absolute values of the  negative principal curvatures of $Y$).  

If $Y$ is non-compact one needs to modify this $f$ by making its second derivative normal to $Y$ to be large as a {\it function} on $Y$.

\vspace {1mm}

{\it Remark.}  The above remains true (and equally obvious) for {\it $k$-mean convex}  hypersurfaces $Y\subset X$, $k=1,2,...,n-1$,  where the traces of the second fundamental forms are positive on  the $k$-planes tangent to $Y$.
\vspace {1mm}

The  notion of $k$-mean convexity  extends from  $C^2$-functions to all continuous ones via {\it linearity} of  the operator $\Delta_k$
by declaring 
a {\it continuous} function $f(x)$  being $k$-mean $\varphi$-convex 
 if it is contained in the {\it localized weak convex hull} of the space of  {\it smooth} $k$-mean  $\varphi$-convex functions. 

In other words $f(x)$ is  $k$-mean  $\varphi$-convex if \vspace {1mm}

{\it $\Delta_k (f)(\tau)- \varphi(\tau)$, understood as a distribution, is representable by a positive 

measure on $Gr_k(X)$.} \vspace {1mm}

Then one defines the set  of  {\it strictly} $k$-mean $\varphi$-convex function as the intersection of the sets of  $(\phi+\varepsilon)$-convex ones, where the intersection is taken over all positive functions $ \varepsilon=\varepsilon(\tau)$ on $Gr_k(X)$.

For example,  a continuous function  $f$ is  strictly $k$-mean convex if there exists a continuous function $\varepsilon(x)>0$
such that  the restriction of $f$ to every  (local) $k$-dimensional submanifold $Y \subset X$ with principal curvatures 
$\kappa_y(Y)$ bounded by $|\kappa_y(Y)|\leq \varepsilon(y)$ is a {\it subharmonic} function on $Y$.

Also one  easily sees that  \vspace {1mm}

\hspace {10mm} {\it  if $f_i$, $i\in I$, are strictly $k$-mean $\phi$-convex functions,  then 
 
\hspace {10mm}   $f(x)=\max_if_i(x)$ is also 
 strictly $k$-mean $\phi$-convex}. \vspace {1mm}

  \vspace {2mm}

 {\it Linearized Definition of  $k$-Mean Curvature Convexity.} A  cooriented hypersurface $Y$  is called 
  {\it strictly $k$-mean $\varphi$-convex }  for $\varphi=\varphi(\tau)$ defined on a neighborhood of the pullback of $Y$ under the tautological map $Gr_k(X)\to X$,  if i$Y$ is representable as the zero set of a continuous strictly $k$-mean  $\varphi$-convex function $f(x)$  defined in some neighbourhood of $Y$.
  
  Here, "cooriented hypersurface" means that  there is a neighbourhood $X_0$ of $Y$ where $Y$ serves as the boundary of a  closed domain $U\subset X_0$. Then our $f$ must be positive inside $U$ and positive outside. We say in this situation  that  $U$  itself {\it strict  
  $k$-mean (curvature) $\phi$-convex (at the boundary).} \vspace {1mm}
  
{\it Mean Curvature Convexity of Functions Revisited. } (Compare   \cite{lawson} \cite {sha})  A continuous function $f(x)$ is called    {\it strictly  $k$-mean  $\varphi$-curvature convex} if, for every point $x\in X$, there exists a convex  $C^2$-function $\beta:\mathbb R\to \mathbb R$ with strictly positive derivative $\beta'>0$, such that the composed function $x\mapsto  \beta\circ f(x)=\beta(f(x)))$  is strictly  $k$-mean  
$\varphi$-convex in some neighbourhood of $x\in X$.

Notice that  
the so defined   strict  $k$-mean curvature  $\varphi$-convexity is 

\hspace {10mm} {\it stable under 
small $C^2$-perturbations of functions.}  

Also, \vspace {1mm}

{\it maxima of families of  strictly  $k$-mean  curvature $\varphi$-convex functions 

are strictly  $k$-mean  curvature $\varphi$-convex,} \vspace {1mm}

\hspace {-6mm} since 
$$\beta\circ\max_if_i=\max_i \beta\circ f_i \mbox { for  monotone increasing $\beta$} $$
 and since  strict  $k$-mean  $\varphi$-convexity is stable under taking maxima.  
 \vspace {1mm}

{\it Remark} Probably, little (essentially nothing?) changes if one allows {\it non-smooth} convex monotone increasing $\beta$
in this definition.
 \vspace {1mm}

{\it $k$-Convexity Lemma.} Let $X$ be a $C^2$-smooth Riemannian manifold and   $U\subset X$  a quasiregular domain with boundary $Y=\partial U$ which is  $k$-mean  $\varphi$-convex} on the  regular locus $reg_Y=Y\setminus sing_Y$ for a positive continuous function $\varphi=\varphi(\tau)>0$. Then \vspace {1mm}

{\it  the minus distance function
$d_U(x)=-dist(x,Y)$ is  $k$-mean    curvature  $(\varphi-\varepsilon)$-convex in the interior of $U$ for some continuous function $\varepsilon=\varepsilon(\tau)$ that vanishes on the pullback of $Y$ in $Gr_k(X)$. 

Moreover,  $\varepsilon(x)$ is bounded  in terms of  $|\inf_x Ricci_x(X)|$ for $x$ running over the ball $B_x(R)\subset X$ for $r=dist(x, Y)$} (For instance,  $\varepsilon(x)=0$ if the Ricci curvature is non-negative in this ball.) 

\vspace {1mm}

The proof is quite simple and,  I guess, is well known in some quarters. Yet,  for the completeness sake, we spell it down in section  5.6. where our argument is essentially the same as that in    \cite  {lawson} and \cite {sha}. \vspace {1mm}
  
 {\it Question.} Is there a meaningful characterization of  Plateau--Stein  $n$-manifolds that
 admit (proper) strictly $(n-1)$-mean convex functions?

 \subsection {Smoothing and  Approximation.}

 Continuous strictly $k$---mean convex function $f$ can be  approximated by   {\it smooth} strictly $k$-mean convex ones,  by convolving with the following  
 
 {\it  Standard $\varepsilon$-Smoothing Kernel.}  Such a kernel is a function   
in two variables, $K_\varepsilon(x_1,x_2)$, $\varepsilon >0$,  on a Riemannian manifold $X$ that is defined with some $\Psi$ by

$$K_\varepsilon(x_1,x_2)  = \lambda(x_1) \Psi(\varepsilon^{-1}dist(x_1,x_2)) , \mbox { }\varepsilon >0,$$  
  for  
  $$\lambda(x_1)=\left (\int_X\Psi(\varepsilon^{-1}dist(x_1,x_2))dx_2\right)^{-1}.$$ 
  
{\it A standard $\varepsilon$-smoothing} of functions on $X$ is 
 $$f(x)\mapsto  f_\varepsilon(x)= \int_X f(x_2)\cdot K_\varepsilon(x ,x_2) dx_2.$$

It is obvious that if $f$ is a {\it continuous} strictly $k$-mean $\phi$-convex function  and $V\subset X$ is a compact subset, then  $f_\varepsilon$ is strictly $k$-mean $\phi$-convex on $V$ for all sufficiently small $\varepsilon>0$. 

It follows that $f$ can be uniformly, and even in the fine $C^0$-topology, approximated by {\it $C^2$-smooth}
strictly $k$-mean $\phi$-convex functions, where, moreover such approximating functions can be chosen equal  $f$ 
on a closed subset $X_0\subset X$ if $f$ itself is  smooth in a neighborhood of $X_0$ in $X$.

Recall that {\it $C^0$-fine approximation} means that the difference between $f$ and an approximating function can be made less than a given strictly  positive continuous function on $X$.

\vspace {1mm}

{\it Curvature Smoothing Corollary.} Let $f(x)$ be a continuous strictly  $k$-mean curvature $\phi$-convex 
function on a $C^2$-smooth Riemannian manifold $X$. Then \vspace {1mm}

{\it $f$ can be $C^0$-finely approximated by smooth strictly  mean curvature $\phi$-convex functions with non-degenerate critical points.} \vspace {1mm}

{\it Proof.} Locally, in a neighborhood $U_x \subset U$  of a given  point  $x$,  such an approximation obtained by finely approximating $  f\circ \beta(x)=f(\beta(x))$, for a suitable $\beta$,  by a smooth $k$-mean $\phi$-convex function on $U_x$, call such an approximation by   $(f\circ \beta)_{appr}$ and then applying the inverse $-\beta$-function, thus approximating $f$ by $f_{appr}=\beta^{-1}\circ (f\circ \beta)_{appr}$. 

Then the global  $C^2$-smooth approximation of $f$ is obtained, by a usual argument  
with a covering of $X$ by  open subsets $U_i$, $i=1,2,...,n+1=dim(X)+1$  where each of  $U_i$ them equals the disjoint union of arbitrarily small subsets. 

Finally,  "$C^2$-smooth" is upgraded to "generic $C^\infty$" by an arbitrarily $C^2$-small perturbation. QED.

 \vspace {1mm}

\section {Splicing, Smoothing and Extending Distance Functions.}

We shall prove in this section the inverse  maximum principles stated  in section

{\it Smoothing Bubbles and Minimal Hypersurfaces.} This smoothing applies, in 

\subsection { Staircase of Distance Functions.} Let $X$ be a Riemannian manifold let   $U_1\supset  U_2 \supset  U_3\supset ...\subset  U_i \supset  ...  \subset X$ be closed domains with boundaries $Y_i=\partial U_i, $ let $d_i(x)$, $x\in U_i$, denote  the  minus distance function from $x$ to the boundary $Y_j=\partial U_j$
and let   some numbers $\delta_i>0$ satisfy
 $$\delta_i>\sup_{x\in Y_{i+1}} dist(x,Y_{i}).$$

\vspace {1mm}

{\it Then  there exists  a negative  proper   continuous function
$h:X \to \mathbb R_-$, such that locally,  in a neighbourhood of every point $x\in X$,  this 
$h$ equals the maximum of the functions $\beta_i\circ d_i$, for $k_x\leq i\leq l_x$, where

$\bullet$  \hspace {1mm}  $l_x$ is the maximal $l$ such that $x\in U_l$;

$\bullet$  \hspace {1mm} $k_x$ is the minimal $k$ such that $ dist (x, Y_k) \leq \delta_k$;

$\bullet$   \hspace {1mm}
$\beta _j$ are smooth monotone increasing  functions, $\beta_k:\mathbb R\to \mathbb R$, with strictly positive derivatives,   $\beta_j' >0$.} \vspace {1mm}

{\it Proof.}  The required max-function $h$ is  determined by its sublevels, call them 
$Y_{-\rho}=h^{-1}(-\infty, -\rho) \subset X$, that  come as  {\it intersections } of  certain  sublevels of the functions $d_i$
that are $(U_i)_{-\rho_i} \subset U_i=d_i^{-1} (\infty,\rho_i]$  for some $\rho_i$ that   must be  continuous strictly increasing functions in $\rho$. 

The essential point is to choose these  $\rho_i$, such that  
  
  {\it if the   boundary of some    $(U_i)_{-\rho_i}$ passes through a point $x$ contained 
  
   in  the boundary of the intersection    $\cap_i(U_i)_{-\rho_i}$, then 
  $dist(x, Y_i)\leq \delta_i$.}

Since   $\delta_i>\sup_{x\in Y_{i+1}} dist(x,Y_{i})$, this inequality can be  obviously  satisfied with some  $\rho_i$  and the proof follows. \vspace {1mm}.

{\it Corollary. Non-Smoothed Inverse Maximal Principle for compact manifolds.} Let $X$ be a compact Riemannian $C^2$-smooth manifold with strictly mean convex boundary.

{\it Then either $X$ contains a compact minimal hypersurface in its interior 

or it admits a continuous negative  strictly mean convex function that 

vanishes on the boundary of $X$.}\vspace {1mm}

{\it Proof.} Shrinking the mean convex  "ends" of $X$ (see 2.4) provides a  finite
 descending sequences of $\phi$-convex bubbles  $U_i$ with a {\it fixed} (albeit very small) {\it strictly positive} $\phi$  and with  arbitrarily small $\sup_{x\in Y_{i+1}} dist(x,Y_{i}) >0$. 
 Then the above $h$ is strictly mean convex being local maximum of distance  
 functions  that are strictly convex by the $k$-convexity lemma
in section 3.3. \vspace {1mm}

{\it Remark. } Bruce Kleiner  explained to me how   a version of this follows by an application   of the mean curvature flow,
but  this   does not seem to be   simpler than our more   pedestrian argument.

\subsection { Proofs of    Inverse Maximum Principles.}

What remans is to justify  step 4  in the proof   of IMP in  section 1.4.

Let, for instance  $X$ be complete $C^2$-smooth Riemannian manifold  that is  connected  and thick
at infinity.   

We already know (see section 2.4)  that if $X$ contains no minimal hypersurface then it 
can be 
exhausted  by compact strictly  mean convex   bubbles $U_j$.

We also know that each $U_j$ can be shrunk via smaller bubbles $U_{ji} \subset U_j$,
$$ U_{j0}=U_j\supset U_{j1} \supset  U_{j2} \supset...\supset  U_{ji}\supset ...,$$
where the minus  distance functions $d_{ij}(x)=-dist(x, \partial  U_{ji})$, $x\in U_{ji}$,
can be "spliced" to continuous  mean  curvature  convex   functions $h_j$ on $U_j$.

If the (positive!) mean curvatures of the  boundary hypersurfaces    $Y_{ji}=\partial U_{ji}$  are bounded from above at all points $x\in Y_{ji}$  by  $\alpha(x)$, where  $\alpha(x)$ is a  (possibly very fast growing)
continuous function on $X$, then, by the usual compactness principle of 
 the geometric measure theory,   some subsequence of $h_j$ converges on all compact subsets in $X$ to the required $h$.

A transparent  way  to achieve the control over $\sup mn.curv(\partial U') $ of a bubble $U'$ inside 
 a given bubble $U$  is to see the construction of $U'$ in terms of   an obstacle (see section 1.4)  that is  a subdomain $V\subset U$ that must be engulfed by $U'$. If 
 the  mean curvatures of the boundary of $V$ at all  boundary points are  bounded by 
$\alpha(x)$, then the same bound will hold for  $\partial (U')$. 

If, for instance, $X$ has Ricci curvature bounded from below, one may take $V=(U_{-\rho})_{+\rho/2}$  where this $V$ (pinched  between 
$U_{-\rho/2}$ and   $U_{-\rho}$) has its mean curvatures  bounded by above roughly by 
$\rho^{-(n-1)}$.

In general, one modifies this by replacing constant $\rho$ by a positive  function $\rho(x)$ on $X$, that must decay, roughly, as   $(1+|R(x)|^{(n-1)})^{-1}$,  for a negative  function    $ R(x)$ that serves as a lower bound for the Ricci curvature of $X$.

 The  curvature of the boundary of such a $V$,  that is obtained by pushing $U$ inward by $\rho(x)$  and then outward by $\rho(x)/2$,  can be easily    bounded by some (very fast growing)  $\alpha(x)$.

This argument, that extends to multi-ended manifolds with the preparations made in sections 2.2-2.4, yields the following non-regular   IMP stated in section 1.4. \vspace {1mm}

{\it Trichotomy Theorem.}  Let $X$ be a complete Riemannian $C^2$-smooth $n$-manifold. Then  (at least) one of the the following three  conditions is satisfied.
\vspace {1mm}

(1) {\it    $X$  admits a proper  strictly  mean curvature convex function

 $h:X\to \mathbb R_+$.}

(2) {\it $X$ 
contains a  
minimal hypersurface $H$  that is  closed in $X$  as a subset 

and  such that    $vol_{n-1}(H)<\infty$.}

(3)  {\it  $X$  admits a  non-proper continuous  strictly  mean curvature 

convex 
function 
$h:X\to \mathbb R$ and also there is a sequence of minimal 

hypersurfaces $H_i\subset X$  
with 
boundaries $\partial H_i$ contained in a fixed compact 

subset
$X_0\subset X$, where 
 these  $H_i$ are   closed in $X$  as  subsets and  such that 
$$vol_{n-1}(H_i)<\infty,  \mbox { \hspace {1mm}   $diam(H_i)\to \infty$ for $i\to \infty$}.$$}

\vspace {1mm}

{\it  Proof of the Regularized Maximum Principles.}  The above functions $h_i$ can be approximated by  $C^2$-Morse functions and minimal   $H$ approximated by slightly concave hypersurfaces according to the smoothing lemma.(see 3.4). This accomplishes the proof of the IMP's  stated in section 1.3.

\vspace {1mm}
  {\it Remark.}  
 It seems, I did not check the details,  the above theorem  remans true with "convex"  replaced everywhere 
 by $\phi$-convex for a given continuous (not even necessarily positive)  function $\phi(x)$,  where the  minimality condition on    $H$ must be replaced by  $mn.curv_x(H)=\phi(x)$, $x\in H$, and where the finiteness requirement for  the $n-1$-volume of $H$ must be replaced by a suitable finiteness condition for  some  $\phi$-area.

 \section {Generalized Convexity.}

We  look  at  the mean convexity from a broader prospective in this section and    we prove the $k$-mean convexity lemma from section 3.4.
All  of what we say is known but dispersed in the literature.

  \vspace {1mm}
  
  \subsection { Smooth and Non-Smooth  Convexity Classes.}
  
   {\it  A coorintation} of a germ of  hypersurface $Y$ at a point  $x $ in a manifold $X$ is  expressed  by calling the closure of  one of the two   "halves"
in the complement   $B_x\setminus Y$, for a small ball at $x$, being  {\it  inside  $Y$}   and the closure of the  other half  {\it outside} $Y$. 
 
Thus,  cooriented  germs at $x\in X$  are   partially ordered.  We agree, thinking of $Y_2$ being {\it more} convex than $Y_1$,  that  

\vspace {1mm}
  
  \hspace {20 mm}$Y_2 \geq_x Y_1$ stands for  $Y_2$  is {\it inside} $Y_1$.

 \vspace {1mm}
Formally, being inside a cooriented
$Y_1$ does not need any coorientation of  $Y_2$. In fact,  $Y_2$  is {\it inside} $Y_1$ implies that 
 $Y_1$  is {\it outside} $Y_2$ {\it only for one} of the two coorientations of $Y_2$. So the above 
 "$Y_2$  is {\it inside} $Y_1$"  tacitly assumes that this does  imply "$Y_1$  is {\it outside} $Y_2$"; moreover,
  if, geometrically, without coorientations,  $Y_1=Y_2$, then  "inside" means that their coorientations are equal as well.

  \vspace {1mm}

  Assume  $X$ is smooth, let $T(X)$ denote the tangent bundle of $X$  and   $\cal H$ be the space of tangent cooriented hyperplanes $H=H_x\subset T_x(X)$, $x\in X$,  that are the tangent spaces to  germs of smooth cooriented hypersurfaces in $X$. Observe that the relation   $Y_2 \geq_x Y_1$
between { \it $C^1$-smooth}  cooriented hypersurfaces implies that they have equal  {\it oriented}  tangent spaces (hyperplanes) at $x$. Accordingly, we may write 
 $Y_2 \geq_H Y_1$ instead of $Y_2 \geq_x Y_1$ for their common cooriented tangent hyperplane $H\subset T_x(X)$.

 \vspace {1mm}

  Given a cooriented hyperplane $H\subset T_x(X)$, denote by  ${\cal Q}_H={\cal Q}_H(X)$ the space of quadratic functions (forms)  $H \to T_x(X)/H=\mathbb R$ and by  ${\cal Q}={\cal Q}(X)$ the space of   ${\cal Q}_H$ over all  $H\subset T_x(X)$, $x\in X$.

   The affine space $aff(  {\cal Q}_H)$ naturally represent  the space of $2$-jets $J_x^2(Y)$ of germs of cooriented hypersurfaces  $Y$ at $x$ that are
tangent to $H$ and  one may speak of the difference 
  $$J_x^2(Y_1) - J_x^2(Y_2) \in  {\cal Q}_H. $$
Obviously,
 $$\mbox {$Y_2 \geq_H Y_1\Rightarrow J_x^2(Y_1) - J_x^2(Y_2)\geq 0$}$$
 where we refer to the usual partial order on the space   ${\cal Q}_H$ regarded as a space of $\mathbb R$-valued functions on $H$. where this implication is reversible for  

\vspace {1mm} 
{\it Strict Order.} The above implication is not, in general, reversible but it is reversible in the strict form:
 $$\mbox {$Y_2 >_H Y_1\Leftrightarrow J_x^2(Y_1) - J_x^2(Y_2)> 0$},$$
where the   {\it strict inequality between germs} signifies that not only   $Y_2 \geq_H Y_1$, but also that this  non-strict inequality is {\it stable} under small $C^2$-perturbations of the germs that remain tangent to $H$.

\vspace {1mm}

If $X$ is endowed with an affine connection, then one may identify $aff(  {\cal Q}_H)$  with $  {\cal Q}_H$; thus, every germ $Y$ is assigned the quadratic form $Q$
on $H=T_x(Y)\subset T_x(X)$ with values in $T_x(X)/H$.  If, moreover,   $X$ is a Riemannian manifold, then 
   $T_x(X)/H$ is {\it canonically} isomorphic to $\mathbb R$, and $Q$ equals the  {\it second   fundamental form of $Y$ at $x$.}

 \vspace {1mm}
  
   Denote by  ${\cal Q}_{aff}(X)$ the space of $2$-jets of cooriented hypersurfaces $Y\subset X$
 and call 
    a subset ${\cal R}   \subset {\cal Q}_{aff}(X)$ a {\it  convexity relation} (of second order) if 
 $$ J_x^2(Y_1)\in {\cal R}\Rightarrow J_x^2(Y_2) \in {\cal R}\mbox { for all germs   } Y_2\geq Y_1,$$
  where $Y_2\geq Y_1$ signifies that both $Y$ have the {\it same} underlying cooriented tangent space   (hyperplane) $H\subset T_x(X)$ where this inequality makes sense.
  
We say that a cooriented   $C^2$-smooth hypersurface $Y\subset X$ 
 {\it satisfies} ${\cal R}$, or it is  {\it ${\cal R}$-convex}, if the $2$-jets of $Y$ are contained in ${\cal R}$ at all points $y\in Y$. 
 
 If $X$ is a Riemannian manifold, then  ${\cal Q}_{aff}(X)= {\cal Q}(X)$ and such a relation is expressed in terms of the second fundamental forms of hypersurfaces.

\vspace {1mm} 

\subsection {  $k$-Convexity  and  $(n-k)$-Mean Convexity.} 

Let $X$ be an $n$-manifold with an affine, e.g. Riemannian, connection and say that a cooriented  $C^2$-hypersurface $Y\subset X$ is $([k_\geq] +[k_>])$-{\it convex} if the second fundamental form of $Y$ with values in $T(X|Y)/T(Y)$, when diagonalized,   has at least  $k_\geq$
nonnegative terms and $k_>$ positive terms. If only one of the two terms in the sum $k=k_\geq +k_>$ is present, one speaks of  {\it $k$-convexity} for $k=k_\geq$ and of {\it strict  $k$-convexity} for $k=k_>$. 

Accordingly, a  domain $V\subset X$ is called $([k_\geq] +[k_>])$-{\it convex} if its boundary is  $([k_\geq] +[k_>])$- convex. 

 For instance, a small $\varepsilon$-neighbourhood of compact smooth submanifold  $P^{n-k-1}$  of codimension $k+1$ in  a Riemannian $X$ is strictly $k$-convex and it is easy to show that every curve-linear subpolyhedron  in $X$ of codimension $k+1$  also admits an arbitrarily small  strictly $k$-convex neighbourhood.

 If $X=\mathbb R^n$,  these  $[k_\geq] +[k_>]$  are {\it the only}    convexity relations that are invariant under {\it affine} transformations of $\mathbb R^n$, where  
  $k=n-1$ corresponds to the ordinary local convexity, while $1$-convex hypersurfaces are nowhere concave.

 The  distinction between  "$\leq$"   and "$<$" is nonessential for {\it compact} $Y\subset \mathbb R^n$,  since, (almost) obviously (see    \textsection 1/2 in \cite {sign})
  
  \vspace {1mm}
every  smooth compact (possibly with a boundary and with a self-intersection)  $k$-convex hypersurface $Y$,
in $\mathbb R^n$ can be $C^2$-approximated by {\it strictly} $k$-convex hypersurfaces  $Y'$ that may be positioned, depending on  what   you wish, inside or outside $Y$. 
 
   \vspace {1mm}
  
{\it Remark/Question.} If $Y$ is non-compact, then a "strict" approximation  of $Y$ by $Y'$  may be possible in one topology, e.g. for $Y'$ being obtained  from $Y$, by a map $f$ with $dist(y, f(y))\leq \varepsilon$ but not in a finer topology where $\varepsilon=\varepsilon(y)\to 0$ for $y\to \infty$.
 Besides, an approximation  of $ ([k_\geq] +[k_>])$-convex hypersurfaces by 
$ ([k'_\geq-l]+ [k'_>+l])$-convex ones  may depend on $l$ and on your positioning $Y'$ inside/outside $Y$.

Is there a comprehensive description of what may happen?

 \vspace {1mm}

Since a generic linear function $f$ on a $k$-convex  domain $V\subset \mathbb R^n$    bounded by $k$-convex hypersurface adds no $l$-handles to  sublevels of $f$  at the critical points of $f$ on $Y$,  

 \vspace {1mm}
{\it a compact $k$-convex  domain  $V\subset \mathbb R^n$    is   diffeotopic to a regular neighbourhood of   $(n-k-1)$-dimensional subpolyhedron  $ P^{n-k-1} \subset \mathbb R^n$.}

 \vspace {1mm}
{\it Questions.} Does there exist  such a diffeotopy $f_t:V\to  \mathbb R^n$,  (that eventually  "shrinks" $V$
to   $ P^{n-k-1}$)   where all intermediate  domains $f_t(V)$, $t>0$ (for $f_0(V)=V$)  are $k$-convex?

 What are  topological possibilities of  $k$-convex domains in the Euclidean $n$-sphere? 
 
Observe that the complement to a disjoint union of  $\varepsilon $-neighbourhoods of two or more equatorial spheres of dimension $k$, $k<n/2$, is strictly $k$-convex; it is contractible to  some $P^{n-1}$  but not to any $P^{n-k-1}$.

  (A more traditional problem concerns  $k$-convex domains $V\subset S^n$, such that, moreover, the complementary domains  $S^n\setminus V$ are $(n-1-k)$-convex.)

 \vspace {1mm}
 
Recall that a  $C^2$-smooth coorineted hypersurface $Y \subset X$  is called  $(n-k)$-{\it mean convex} if the traces of the second fundamental form of $Y$ restricted to the  tangent $(n-k)$-planes $H^{n-k} \subset T(Y)$ are non-negative.  In other words,
 
{\it the principal curvatures of $Y\subset X$, say $\kappa_1\leq  \kappa_2\leq... \leq\kappa_{n-1}$ satisfy
$$\kappa_1+  \kappa_2+...+ \kappa_{n-k}\geq0\mbox { at all points $y\in Y$.}$$}
 (If $k=1$ this means that the  $Y$ is convex and if $k=n-1$ this says that the mean curvature of $Y$ is non-negative.)

Accordingly,  {\it  strict $(n-k)$-mean convexity}  requires  this inequality to be strict, i.e. all  traces to be positive. 

{\it Question.} Can every  closed  $(n-k)$-mean convex hypersurface in $\mathbb R^n$  be approximated by  {\it strictly} $(n-k)$-mean convex ones?  

(This is easy for $k=1$ and $k= n-1$,  but I see no immediate proof it for other $k$.  Am I missing something obvious?\footnote{Bruce Kleiner pointed out to me that such approximation is possible with the mean curvature flow. }

 \vspace {1mm}
 
 Clearly, (strictly) $(n-k)$-mean convex hypersurfaces are (strictly) $k$-convex, and

 {\it every embedded closed  $k$-convex hypersurface  in the Euclidean space  $\mathbb R^n$  is isotopic to a strictly $(n-k)$-mean convex one} (since it can be brought to a neighbourhood of $k$-subpolyhedron $P^k \mathbb R^n$).

 But this 
 if far  from being true, even on the homotopy level, in non-Eulidean spaces.
 
  \vspace {1mm}

{\it Mean Convex Surgery.}  Let $V\subset X$ be a smooth $(n-k)$-mean convex domain and let 
$B^l\subset X$ be a smooth disk that all lies outside $V$ except for its    boundary sphere $S^{l-1}=\partial B^l \subset Y=\partial V$, where we assume  (just for the  civility sake) that $B^l$ meets $Y$ {\it normally} i.e. under the angle$=\pi/2$ along 
  $S^{l-1}= B^l\cap Y$.

Let us slightly thicken $B^l$ by taking its $\varepsilon$-neighbourhood, denoted $_\varepsilon B^l\subset X$, and observe, assuming $\varepsilon>0$ is sufficiently small, that 

$\bullet$ the union
$V\cup  {_\varepsilon B}^l$ has smooth boundary except for a   $\sim\pi/2$ corner along the boundary of a small tubular neighbourhood of $S^{l-1}\subset Y$; 

$\bullet$ the new  smooth part of the boundary of  $V\cup  {_\varepsilon B}^l$, that is
$$ \partial (V\cup  {_\varepsilon B}^l) \setminus \partial V=(\partial( _\varepsilon B^l)\cap (X\setminus V)$$ 
 is $(l+1)$-mean convex.

\vspace {1mm}

{\it If $l \leq  n-k-1$, then the corner in $V\cup  {_\varepsilon B}^l$ can be $(n-k)$-mean convexly smoothed.} 
 
\vspace {1mm}

{\it About the  Proof.} The boundary of the union  $V\cup  {_\varepsilon B}^l$ is concave along the corner and the obvious smoothing of $V\cup  {_\varepsilon B}^l$ does not give us a $(n-k)$-mean convex domain. However the $(n-k)$-mean  curvature of the  boundary of the $\varepsilon$-tube around $B^l$ for   $l \leq  n-k-1$ tends to $+\infty$  for  $\varepsilon \to 0$.   This  "infinite excess of positivity"  allows one    to   to construct strictly  $(n-k)$-mean  smoothing similarly but easier  than how it was  done in \cite {GL} for scalar curvature.

\subsection { Convergence Stability.} The limit behavior of  {\it embedded}   ${\cal R}$-convex hypersurfaces is opposite to what is demanded   by the {\it $C^0$-dense $h$-principle}: the spaces of such hypersurfaces are {\it closed} rather than dense in the $C^0$-topology for closed subsets ${\cal R} \subset \cal Q$.

Moreover, let ${\cal R}  \subset { \cal} Q(X) $ be a  closed convexity relation and let 
  $Y\subset X$ be a  $C^2$-smooth  cooriented hypersurface that is  closed in $X$ as a subset. Let $U_i\supset Y$, $i=1,2...$, be a sequence of neighbourhoods such that 
  $\bigcap_i U_i=Y$ and let $Y_i\subset U_i$ be smooth cooriented  hypersurfaces closed in $U_i$ as subsets,  the closures of which do not intersect the boundaries of $U_i$ and  that 
  separate the  components of the boundaries $\partial U_i$ in the same way $Y$ does. In other words, $Y_i$ are homologous to $Y$ in $U_i$  (in the sense of homology with infinite supports if $Y$ is non-compact.)
  \vspace {1mm}
  
{\it If all   $Y_i$ satisfy a closed convexity relation $\cal R$ then $Y$ also satisfies $\cal R$.}
   
   \vspace {1mm} 
  
 {\it Proof.} In fact let $Q_0$ be the jet of $Y$ at some point $Y_0$ and $\Omega_0\subset \cal Q$
  be an   arbitrarily  small neighbourhood of $Q_0$. Then, by the {\it weak} (and obvious)
  {\it maximal principle},  every $Y_i$ for all  $i\geq i_0=i_0(\Omega_0)$  contains a point $y_i^-$
 such that the  $2$-jet  $J^2_{y_i^+}(Y_{i})\in  \cal Q$  (or a germ at this point, if you wish) satisfies 
 $$ J^2_{y_i^-}(Y_{i}) \leq \omega_i \mbox { for some } \omega_i\in \Omega_0$$
  (as well as a point  $y_i^+$, where $$ J^2_{y_i^+}(Y_{i}) \geq \omega'_i \mbox { for some } \omega'_i\in \Omega_0.)$$
Q.E.D

\vspace {1mm} {\it (Counter)examples.}   (a)  Every curve $Y$ in the plane can be (obviously) $C^0$-approximated by  locally convex {\it immersed} curves $Y_i$.    By the above, these $Y_i$ must have lots of double points as $Y_i$ come close to the region where $Y$ is concave.   

(b) Similarly, according to Lawson and Michelson \cite {lawson1} every co-orientable immersion $f$ of an $(n-1)$-manifold  $Y$ to a Riemannian $n$-manifold $X$, can be  $C^0$-approximated by immersions $f_i: Y\to X$ with positive mean curvatures.

(The building blocks of $f_i$ are    finite coverings maps onto  the boundaries  of $\varepsilon$-neighbourhoods  of $(n-2)$-submanifolds in $X$, where, observe, such boundaries  have mean curvatures $\sim \varepsilon^{-1}$ for 
$\varepsilon \to 0$.)
\vspace {1mm}

(c) In contrast with the above, if $k >n/2$, then\vspace {1mm}

 {\it every  closed cooriented strictly $(n-k)$-mean convex hypersurface $Y$ in  a complete Riemannian manifolds $X$ with non-negative sectional curvatures 
bounds a compact Riemannian manifold $U$, i.e. $\partial U=Y$, such that the immersion $Y\to X$ extends to an isometric immersion $U\to X$.}

(This $U$ is contractible to its $k$-skeleton, since the minus distance  function $u\mapsto -dist_U(u,\partial U=Y)$  admits an approximation by an $(n-k)$-mean convex Morse function on $U$  that provides an isotopy of $U$ to a regular neighbourhood of a  $k$-dimensional subpolyhedron $P^k\subset U$ \cite {sha}.)

On the other hand, $k$-convex  hypersurfaces for $k<n-1$  in  general  {\it non-flat}  $n$-manifolds $X$ with non-negative  curvatures {\it do not  necessarily bound} immersed $n$-manifolds in $X$.  But this is true in the presence of many  "movable" totally geodesic submanifolds in $X$ by the  the Euclidean argument from   \S 1/2 in \cite {sign},  where  the standard  examples of such manifolds    are  Riemannian cylinders  $X=X_0\times \mathbf R$ and complete 
 simply connected  $n$-spaces  $X$ of {\it constant negative} curvature.

\vspace {1mm}  
   
 {\it Questions. } What are  possible  topologies  of  (embedded and immersed) $k$-convex hypersurfaces in the Euclidean  $n$-sphere?
 
  Are there any constrains on immersed  $(n-k)$-mean convex hypersurfaces in  the Euclidean $n$-space for $n\geq 2k$?

 \vspace {1mm}  
 
 The  convergence stability suggests  that that the notion  of $\cal R$-convexity can be   extended to non-smooth subsets.  The cheapest way to produce non-smooth examples   starting with the class ${\cal U}$   of smooth {\it $\cal R$-convex  domains $U \subset X$},   i.e. having smooth   $\cal R$-convex boundaries, is  enlarging/completing   ${\cal U}$ by some/all of  the following three operations over subsets.
 
 $ [\uppitchfork]_{<\infty} $  Locally finite intersections of smooth domains  $U_i\in \cal U$ with mutually transversal boundaries. 
 
 $ [\bigcap\downarrow]_\infty $ Intersections of infinite  {\it decreasing} families of subsets.

  $[\bigcap]_{\leq \infty}$   Finite and infinite  intersections of smooth domains, that is essentially the same as  $[\uppitchfork]_{<\infty}+ [\bigcap\downarrow]_\infty $.
 
  $ [\bigcup\uparrow]_{-\infty} $ Union of infinite  {\it increasing} families.

 \vspace {1mm}
 
 {\it Question.} Given, say an open,  convexity relation ${\cal R}\subset {\cal Q}(X)$, let ${\cal C}({\cal R})$ denotes the class of all compact subsets in $X$ obtained from compact  smooth $\cal R$-convex domains $U\subset X$ by some of the above  operations, e.g. by   $[\bigcap]_{\leq \infty}$, i.e. by taking infinite   intersections of  compact  smooth $\cal R$-convex domains $U$.

 Is there any, not necessarily exhaustive, characterization of subsets in    ${\cal C}({\cal R})$
 in terms of  ${\cal R}$?
   
Is every  ${\cal R}$  {\it  uniquely} determined by ${\cal C}({\cal R})$?  
   
For instance, which  {\it Cantor    sets} $C \subset \mathbb R^n$ are representable as infinite intersections of disjoint  finite unions of compact {\it convex} subsets?
   
  Clearly, this is possible if the Hausdorff  dimension of $C$ satisfies $dim_{Hau}(C)<1$,
  but "generic" subsets $C$ with $dim_{Hau}(C)>1$  admit no such representation.
   
  In fact, the geometry of a  Cantor set  $C\subset \mathbb R^n$ at a point $x\in C$ may be characterized by 
  the minimal possible "oscillatory complexity", $osc_\varepsilon=osc_\varepsilon(C,x)$, $\varepsilon >0$, e.g the {\it total curvature} $curv_\varepsilon$  (that is the ${n-1}$-volume of the tangential Gauss map counted with multiplicity) of  the boundaries of smooth  neighborhoods 
  $U_\varepsilon\subset \mathbb R^n$ of $x$ such that  $diam(U_\varepsilon)\leq \varepsilon$ and where the boundaries  $\partial U_\varepsilon$ do not intersect $C$. 

 It seem "most" Cantor   sets in $\mathbb R^n$, $n\geq 2$,    (I   checked this only for a few particular classes of sets)  have  
$osc_\varepsilon\to \infty $, e.g. 
 $curv_\varepsilon \to \infty $   for $\varepsilon \to 0)$ and they do not belong to any 
 convexity class  ${\cal R}  \subset {\cal Q}$, unless  ${\cal R} $ equals
 $\cal Q$ minus a  "very thin"   subset).

 \vspace {1mm}  
   
{\it Convergence Stability for $k$-Mean Convexity for Functions.} Since this convexity is defined by {\it linear} inequalities on the (second) derivatives of functions $f$ it is stable under all kinds of weak limits
 and  it non ambiguously extends to continuous functions as we saw in section 3.3. 

\subsection {   Riemannian Curvature  Digression.} The above is a baby version of the following Riemannian problems.
   
   Given two  $2$-jets, or germs $g_1$ and $g_2$ of Riemannian metrics at a point $x$ in a smooth manifold $X$, write $g_1 \preccurlyeq g_2$, if the two have equal $1$-jets and their sectional curvatures satisfy

   $$curv_\tau(g_1)\geq curv_\tau(g_2)\mbox { for all tangent  $2$-planes }  \tau\subset T_x(X).$$
   
For example,  metrics  with "large amount" of positive curvature  are regarded as   small.   
   
 A {\it  lower  curvature relation/bound} ${\cal B}$ is a subset of $2$-jets $g$ of Riemannian metrics  at the origin in $\mathbb R^n$, such that \vspace {1mm}

$\bullet$  $1$-jets of $g$  equal the $1$-jet of   the Euclidean metric; 
 
$\bullet$   if    $ g_2\in {\cal B}$ and   $g_1 \preccurlyeq g_2$ then   $  g_1\in {\cal B} $;
   
 $\bullet$  the subset $\cal B$ in the space of jets is invariant under orthogonal 
 
   transformations of $\mathbb R^n$.  \vspace{1mm}
   
  The latter condition allows one to invariantly speak of {$\cal B$-positive metrics} on {\it all } smooth $n$-manifolds $X$ that are, in other words,  {\it  Riemannian $n$ manifolds that satisfy $\cal B$} (compare with \S7 in \cite {sign}).
 
   The fundamental questions are as follows.  \vspace{1mm}
  
  [A]  Given  $\cal B$    what is the weakest topology/convergence  ${\cal T}={\cal T}({\cal B})$  in the space of Riemannian manifolds,  such that the limits of  $\cal B$-positive manifolds are  $\cal B$-positive? 
   
[B] What are  singular  $\cal B$-positive metrics  spaces?
  
[C]  What are    $\cal B$ for which the above two questions have satisfactory answers? \vspace{1mm}

   If $\cal B$ consists of the metrics with a given bound {\it all } sectional  curvatures, then 
    the (best known)  answer to  [A] is  the {\it Hausdorff convergence} of metric spaces 
    and [B] is essentially resolved  by the theory of Alexandrov spaces.
   
   The starting point of   the theory for spaces with a lower bound on the {\it Ricci}  curvature is the (almost obvious)  stability of the  inequality $Ricci(g)\geq const \cdot g$ under $C^0$-limits of Riemannian metrics on a given  underlying (and unchangeable)  smooth manifold $X$ while the general  theory, albeit not fully established,  is 
   well underway, see \cite {che}, \cite {vil}, \cite {olliv} and references therein.
   
 The most tantalizing relation $\cal B$  is expressed with the {\it scalar curvature} by $scal(g)\geq cost $   where even the $C^0$-limit  stability is not fully established  and where  some possibilities are suggested by the {\it intrinsic flat distance} \cite {wenger2}.

Nothing seems to be known about other $\cal B$, e.g. those  encoding some  positivity
of the curvature operator, e.g. positivity of the  {\it complexified sectional curvature}, see \cite {m-m}, and \S7 in \cite {sign}.

 \subsection {Cut Locus, Focality and Curvature Blow-up.}  
 Let us  see what happens to convexity under equidistant deformations of a hypersurface  $Y\subset X$, where an attention must be paid to 
  singularities on the {\it cut locus}  that may be aggravated by the presence of  {\it focal points}.\vspace {1mm}

Recall that the  cut locus $cut(U)\subset X$ of  an open subset $U\subset X$ (or of a closed domain $U$) with respect to $Y=\partial U$ is defined as  the closure of the set of points   $ u\in U$ that have more than one  {\it ancestor} in $Y$, where a point $x$ in the closure of $ U$  is called a {\it $d$-ancestor}, for $d= dist(x,u)$, or just "ancestor"  of a point  $u \in U$, with $u$ being called a {\it $d$-descender}, or "descender" of $x$, if $  dist(x,u)=dist(u, Y)- dist (x,Y).$

 Assume that $X$ is a complete  $C^2$-smooth Riemannian manifold and recall a few obvious properties of $cut(U)$.
 
 \vspace {1mm}

 If $Y=\partial U$ is a $C^2$-hypersurface, than the cut locus of $U$ does   not intersect $Y$ and the $\rho$-equidistant hypersurfaces denoted  $Y_{-\rho} \subset U$ are $C^2$-smooth away from $cut(U)$,
 i.e. the complements $Y_{-\rho}\setminus cut(U)$ are $C^2$-smooth (locally closed) hypersurfaces in $U$.

 \vspace {1mm}

If   the boundary  $Y=\partial U$ is 
   {\it $C^1$-smooth}, then 
   an  $x \in Y$ is  $\rho$-accessible from $U$ if and only  if
  the  geodesic segment of  length $\rho$   issuing from $x$   normally to $Y$ inward $U$  either {\it does not} intersect  $cut(U)$, or, if it meets  $cut(U)$, then  only at its terminal in $U$. 
 
  \vspace {1mm}

 All open $U\subset X$ satisfy   (by  a simple  {\it $\check{C}ech$   homology }  argument).
  
  \vspace {1mm}
  
  \hspace {4mm}    $ \rho \leq dist(y,cut(U)) \Rightarrow$ $y$ {\it is $\rho$-accessible from $U$  for all $y\in Y=\partial U$.}

\vspace {1mm}
  
Consequently, if $U$ is $C^2$-quasi-regular, then 
  $Y\setminus cut(U)$ is $C^2$-smooth.

\vspace {1mm}

 {\it Focal Points.}   Let $y_0\in Y$ be an  ancestor of $u_0\in U$, i.e. a {\it (global) minimum} point of   
 the function $y\mapsto dist(y,u_0)$ on $Y$. Assume $X$ is complete and let 
 $\gamma=\gamma(s)$ in $X$ be a geodesic ray  issuing from $x_0$ inward $U$, such that 
$$\gamma(s_0)=u_0 \mbox { for  } s_0=dist (u_0, y_0),$$
where $s\geq 0$ denotes the  geodesic length parameter. 
(If $Y$ is { \it $C^1$-smooth hypersurface}  at  $y_0$  then  $\gamma$ is unique being normal to $Y$.)  

The point $u_0$ is called {\it non-focal} for $y_0$ along $\gamma$ if $y_0$ remains a {\it local} minimum of the function   $y\mapsto dist(y,u)$ on $Y$  as we slightly move   along $\gamma$ inward,  i.e.
for $u= \gamma(s_0+\varepsilon)$ and  all sufficiently small $\varepsilon>0$.

 \vspace {1mm}

In other words, the   $(s_0+\varepsilon)$-ball in $X$ around $u_\varepsilon\in U$, say 
$$B_{u_\varepsilon} (  s_0+\varepsilon)\supset B_{u_0}(s_0) \subset U,$$
 is  "contained in $U$   at $y_0$", i.e. the intersection of $B_{u_\varepsilon} (  s_0+\varepsilon)$ with a small neighbourhood of $y_0$ in $X$ is contained  in $U$.

 \vspace {1mm}

Notice that focal/non-focal for $y_0\in Y$ depends only on  the geometry of $Y$ in a small  neighbourhood of $Y_0$ plus  on how one defines "inward".  Thus, one can extend the above definition by taking an arbitrarily small neighbourhood $B_0 \subset X$ of $y_0$, (e.g. a small $\varepsilon$-ball  around $y_0$), letting 
$$U_0=X\setminus (B_0\bigcap (X\setminus U))\supset U$$  
and defining focal/non-focal along geodesic segments in $U_0$ that starts at $y_0$ and may go beyond $U$.

\vspace {1mm}

If $Y$ is a  {\it $C^1$-smooth} hypersurface  and $y$  is an ancestor of $u$ with $dist(u,y)=dist (u,Y)=\rho$  then the $\rho$-sphere around $u$, say $S_u(\rho)=\partial B_u(\rho) \subset U$, that contains $y$ 
is $C^2$-smooth at $y$, provided our Riemannian metric is $C^2$-smooth.
If, moreover, $Y$ is a  {\it $C^2$-smooth} hypersurface, then the second  fundamental form $Q_Y$ of  $Y$ at $y$ is  minorized  by the form  $Q_S$  at $y$, i.e.  $Q_Y-Q_S$ is  {\it negative semidefinite}
 since $ B_u(\rho) \subset U$. (Our sign convention for $Q's$ is the one for which the  boundaries  of convex subsets $U \subset X$
  have {\it positive} definite second fundamental forms $Q$.)

Obviously, $u\in U$ is {\it non-focal} for $Y=\partial U$, (along the minimal  geodesic segment between the two points) if and only if the quadratic   form  $Q_Y-Q_S$ is  {\it negative definite}.

 Denote by  $foc(U)\subset U\cap cut(U)$ the subset of {\it the focal points} where $u$ is called focal 
 if it is focal for some  ancestor of $u$ in  $Y=\partial U$ and observe that

\vspace {1mm}

{\it  If $Y$ is {\it $C^2$-quasi-regular}, e.g.  $C^2$-smooth,  then the  subset $foc(U) \subset U$,   is closed in $U$.} 
(This is not, in general, true for $C^1$-hypersurfaces $Y.$)  \vspace {1mm}

 \vspace {1mm}
 
 The appearance of focal points can be  seen in terms of the hypersurfaces $Y_{-\rho} \subset U$ 
 equidistant to  $Y=\partial U$ as follows. Join a point $u_0\in Y_{-\rho}$
 with one of its ancestors, say $y_0\in Y$ by a minimal geodesic segment $\gamma$ in    the closure of $U$, where $length(\gamma)= \rho$,
 and observe that the hypersurfaces $Y_{-\rho+\varepsilon}$, $0<\varepsilon\leq \rho$, are $C^2$-smooth at the points  $u_{+\varepsilon}=Y_{-\rho+\varepsilon}\cap \gamma$, provided $Y$ is $C^2$-smooth at $y_0$. Then, 
  
  \vspace {1mm}
 
 {\it the second fundamental forms $Q_{\varepsilon}$ of  $Y_{-\rho+\varepsilon}$ at the points $u_{+\varepsilon}$
 are uniformly 
 
 bounded from below.}
 
  \vspace {1mm}
 
  \hspace {-7 mm} If  $u_0$ is a {\it non-focal} for $y_0$ then the forms $Q_\varepsilon$ are also bounded from {\it above}; moreover,  the  hypersurfaces $Y_{-\rho+\varepsilon}$ can be locally, around $\gamma$, be included into a $C^2$-smooth family of local equidistant hypersurfaces to a small neighbourhood of $y_0 \in Y$. 

But if    $u_0$ is {\it focal for $y_0$} then these forms "blow up" for $\varepsilon \to 0$ as follows.

The ($(n-1)$-dimensional) spaces $ T(\varepsilon)$ normal  to $\gamma$ at the points $u_{+\varepsilon}$, that serve as tangent spaces to $Y$ for $\varepsilon >0$, admit orthogonal splittings   $T(\varepsilon)= T_0(\varepsilon)\oplus T_1(\varepsilon)$, where these $T_0(\varepsilon)$ and $T_1(\varepsilon)$  continuously depend  on $\varepsilon \in [0,\rho]$ and such that

$\bullet $ the subspace   $T_0(\varepsilon=\rho) \subset T_{y_0}(Y)$  equals the kernel of the above difference form $Q_Y-Q_S$ at $y_0$;

$\bullet $  the forms $Q_\varepsilon$ restricted to $T_1(\varepsilon)$ are continuous for  all  $0< \varepsilon\leq \rho$ and they continuously extend to the space $T_1(\varepsilon=0)$. 
 
$\bullet $ The forms  $Q_\varepsilon$  on the subspaces $T_0(\varepsilon)$ tend to $+\infty$ for $\varepsilon \to 0$. In fact, the values of $Q_\varepsilon$ on the unit vectors in $T_0(\varepsilon)$ is of order $1/\varepsilon$.

 \subsection {$C^2$-Approximation with Corners.}

  We show here how equidistant hypersurfaces to a quasiregular $Y$ can be approximated by piecewise smooth hypersurfaces with one sided  controls on their curvatures.\vspace {1mm}

 Let $U\subset X$ be a $C^2$-quasiregilar open subset (domain) with boundary   $Y=\partial U$  in a complete Riemannian manifold $X$, let   
   $Y_{-\rho}=\partial  U_{-\rho}= U_{-\rho}^\leq \subset U $, $\rho>0$,  be     the   equidistant hypersurface, where, as earlier,    
  $U_{-\rho}^\leq $ denotes    the set of $u\in U$ where $ dist(u,Y)\leq \rho$. \vspace {1mm}

{\it $\smile$Approximation Lemma.} Given   $\varepsilon>0$ and $0<\rho'<\rho$, there exists a domain $U'= U^{\smile\varepsilon}_{\rho'}$ in $X$ such that  
$$  U_{-\rho} \supset   U'\supset U_{-(\rho+\varepsilon)},$$
and such that 

\hspace {14mm}{\it the boundary $Y'=\partial U'$  is  piecewise $C^2$-smooth}. 

 In fact,
 there are $C^2$-diffeomorphisms $D_i:X\to X$  such that 
 
 {\it $D_i(U_{-\rho})$ do not intersect the singular locus of $Y_{-\rho'}$  and 

$Y'$  equals the union of the $D_i$-pullbacks of $Y_{-\rho'}$ 
$$Y'=\cup_i D_i^{-1}(Y_{-\rho'}).$$}
Moreover, if $\rho-\rho'$ is small, then these 
 $D_i$ are $C^2$-close to the identity map $X\to X$;
consequently, 

\hspace {10mm} {\it the curvatures of the smooth pieces  $D_i^{-1}(Y_{-\rho'})$  are close 

\hspace {10mm} to the curvatures of their $D_i$ images in $Y_{-\rho'}$.}

  \vspace {1mm}
 
 {\it Proof.}   Let $\gamma \subset X$ be a minimal geodesic segment. Then, obviously, there exists   smooth vector field 
 $V_\gamma(x)$ on $X$ that is tangent to $\gamma$  where it equals the unit field 
   directed from $x_0 $ to $x_1$ and such that the norm of $V_\gamma$ satisfies 
   $$\mbox{ $||V_\gamma(x)||<1$ for  $x\notin \gamma$ and  $\limsup_{x\to \infty} ||V_\gamma(x)|| <1$}.$$
   Integrate $V_\gamma$  for the flow time   equal $length(\gamma)$;  thus, obtain 
  a $C^2$-diffeomorphism  $D_\gamma:X\to X$  such that

   \hspace {2mm}   {\it $D_\gamma$ sends one end of $\gamma$, say $x_0$, to the other one, called $x_1=D_\gamma(x_0)$,} 
    
 \hspace {-6mm}   where this diffeomorphism is {\it sharp at $\gamma$} in the sense that  \vspace {1mm}
    
   \hspace {20mm}   $dist(x,D_\gamma(x))<length (\gamma)$  for all 
  $x\nin \gamma$,\vspace {1mm}
    
 \hspace {-6mm}  and where one can achieve  a map     
    $\gamma\mapsto D_\gamma$   to be  {\it continuous for the $C^2$-topology}  in the space  of diffeomorphisms.

   \vspace {1mm}

  {\it Remark.}
It is easy to arrange the maps $D_{\gamma} :X\to X$ 
such that their differentials  $T_{x_0}(X)\to T_{x_1}(X)$ are {\it isometries } for all $\gamma$.
Moreover, if $X$ has {\it positive} sectional curvatures,
one can make  $D_\gamma$  {\it second order isometries} at these points, i.e. such that every geodesic through $x_0$ goes to a curve with zero curvature at $x_1$. However, this is impossible for manifolds of {\it negative} curvature.

 \vspace {1mm}

Now, let $\delta>0$  be  very small (depending, in particular, on $\varepsilon$),  take
all minimal segments $\gamma$ between the points $y\in Y_{-\rho}$ and its $(\rho'+\delta)$-anscestors in $U$
and let 
$$U'= U_{-\rho}\setminus \cup_\gamma D_\gamma^{-1}(X\setminus U_{-\rho'}).$$  
 Finally, take a sufficiently dense  locally finite set of geodesic segments, say $\{\gamma_i\}$,
 and take $D_{\gamma_i}$ for the required diffeomorphisms $ D_i$.

 \vspace {1mm}
This $\smile$-approximation implies in particular that  the distance function $d$ to the boundary $Y$ of $U$ can be approximated by  the maximum of  {\it smooth} distance functions with their second partial derivatives close to those of $h$ at nearby points. It follows that all  $k$-convexity bounds extends from smooth to non-smooth points of $d$. In particular the
$k$-convexity lemma  from section follows from this  $\smile$-approximation since one, obviously,  has a uniform bound on the "bending"  $\beta$ in this case. 

\vspace {1mm}

{\it On External Approximation.} The above piecewise smooth hypersurfaces $Y_{-\rho}^{\delta\smile \Delta}$ that approximate the boundary $\partial U^<_{-\rho} $
 are positioned {\it inside}  $U^<_{-\rho}$. Probably, there is no similar approximation by hypersurfaces lying {\it outside} but this is obviously possible if $Y_{-\rho}$  is compact: just apply the inside approximation to $Y_{-\rho'}$ for  $\rho'<\rho$ and let  $\rho'\to\rho$.

\vspace {1mm}

 \subsection{Cornered Domains and Smoothing the Corners.}

\vspace {1mm}

Let us indicate here a geometric alternative to the smoothing operators we used in section 3.4.

     {\it A cornered domain} of class $C^k$ in  a   $C^k$-smooth $n$-manifold  $ X$   is a closed subset $V$, such that every {\it boundary} point $v$ in $ V$ admits a neighbourhood $U(v)$ in $V$ that is $C^k$-diffeomorphic to the intersection of $k\leq n$ mutually orthogonal halfspaces in $\mathbb R^n$. 
  
 {\it  The regular part} of the boundary of $V$, denoted    $reg_{\partial V}\subset \partial V$, consists of those $v$, where $U(v)\subset V$ is diffeomorphic to a half space, i.e. $k=1$.
  
{\it The   $(n-1)$-faces $W_i$, $i\in I$, of $V$} are  the closures of the connected components of   $reg_{\partial V}\subset V$
where, obviously, 
 $$\bigcup_{i\in I}W_i=\partial V.$$  
(Sometimes, one takes finite unions of disjoint  connected  components for faces.)

 {{\it Corners or $(n-2)$-faces of $V$} are,  by definition,  non-empty  pairwise intersections of  $(n-1)$-faces, 
 $$W_{i_1i_2}=W_{i_1}\cap W_{i_2},$$
 (Since we assume the corner structure being "simple", there is no non-empty intersection $W_{i_1}\cap W_{i_2}$  of dimension $< n-2$. On the other hand,  corners may be disconnected.)

It is easy to see that a cornered  $V$  
equals an intersection
 $$V=\bigcap_{j\in J}V_j,$$
where 

\vspace {1mm}

(1) $V_j\subset X$ are  $n$-submanifolds  with {\it smooth} boundaries $\partial V_j,$

(2 ) all intersections between
 ($k$-tuples of) $\partial V_j$ are {\it transversal},

(3 ) there are at most {\it finitely many} boundaries that intersect a given {\it compact} subset in $X$.
   
\vspace {1mm}
These (1)-(3)  imply that
the intersections of $\partial V_j$ with $V$ equal finite union of disjoint $(n-1)$-faces of $V$ where these $V\cap  \partial V_j\subset   \partial V_j$ are cornered domains in  $\partial V_i$.
If one wishes,  one may let $J=I$ and choose  $V_i$ such that $V\cap \partial V_i= W_i$.

 \vspace {1mm}

If $X$ is a Riemannian manifold than one may speak of the dihedral angles between pairs of 
  $(n-1)$-faces along $(n-2)$-faces. Clearly all these angles $\angle(W_{i_1},W_{i_2})$, that are {\it continuous  functions} on  
 $W_{i_1i_2} =W_{i_1}\cap W_{i_2}$, are bounded by  $$\angle(W_{i_1},W_{i_2})<\pi.$$
 
{\it Essential Example.}  A generic  $C^\infty$-perturbation of the smooth pieces of  $U'$ that approximate$U_{-\rho}$ in the previous section  turn $U'$ into a cornered domain.

 \vspace {1mm}
  
 Since
 the corners of  $V$ are convex  for the dihedral angles $<\pi$ one expects 
 that cornered domains  $V$  admit  approximations by  {\it smooth} domains 
 that are, up to an arbitrary small error, "as convex" as the faces $W_i$ of $V$. Indeed,  this is possible for quite a few, classes of {\it convexity relations} (see next section)  including  strict $k$-mean convexity where the picture is most transparent for the mean curvature convexity. 
 
  \vspace {1mm} 
  
  {\it Corner  Smoothing Lemma.} Let $X$ be a $C^2$-smooth Riemannian $n$-manifold, $\phi:X\to \mathbb R$ a continuous function and $V$ a cornered domain of class $C^2$ such that the  mean curvatures of the regular part of the boundary $ \partial V\subset V\subset X$   satisfy, 
 $$mn.curv_x (\partial V)>\phi(x)\mbox{ for all }x \in reg_{\partial V} \subset \partial V.$$
 
 {\it Then, for an arbitrary neighbourhood $\Delta\subset V$ of the boundary $\partial V\subset V$, 
  there exists a domain $V'\subset V$ with $C^2$-smooth boundary, such that
  $$\mbox { $\partial V'\subset \Delta$  and  } mn.curv_x(\partial V')>\phi(x)\mbox{ for all }x \in\partial V';$$
moreover, the normal projection   $\partial V \to \partial V'$ is a $C^2$-diffeomorphism on every  $(n-1)$-face of $V$.}

  \vspace {1mm}
{\it  Proof.} Let $\partial V$ be compact and let $\partial _{+\delta} V\subset X$ be the boundary of the $\delta$-neighbourhood of $V$ for a small $\delta>0$. Clearly, $\partial _{+\delta} V$ 
 is a $C^1$-smooth hypersurface, that is, moreover, piecewise $C^2$. 
 
 This is seen with the normal projection  $\partial _{+\delta} V\to \partial V$ that sends every $C^2$-piece of    $\partial _{+\delta} V$ onto an $m$-face of $V$, for some $m=1,2,...,n-1$. If $m= n-1$,
 then the mean curvature of this piece is $\delta$-close to that of the corresponding $(n-1)$-face
 and if $m< n-1$ then the mean curvature is $\sim \delta^{-1}$.
 
  Thus,  the mean curvature of $\partial _{+\delta} V$ is a piecewise continuous function on    $\partial _{+\delta} V$   that satisfies
 $$mn.curv_x(\partial _{+\delta} V)> \phi(x)\mbox{  for all sufficiently small } \delta>0  \mbox{ and  all } x\in  \partial _{+\delta} V.$$
Now, observe that  $\partial _{+\delta} V$ equals the $\delta$-level of the distance function
$d(x)=dist(x, V)$ and let $d_\varepsilon(x)$ be the average of $d(x)$ over the $\varepsilon$-ball  $B_x(\varepsilon)\subset X$ for a small $\varepsilon>0$. 

Since the second differential of $d(x)$ is a bounded measurable function and $||grad(d)||=1$
the $\delta$-level say $\partial _{+\delta, \varepsilon} V \subset X$ of $d_\varepsilon$ is a $C^2$-smooth hypersurface that $C^1$-converges to $\partial _{+\delta}$ for $\varepsilon  \to 0$.
Since the mean curvature of a  level  of a function is {\it linear} in the second derivatives of the function
 the mean curvatures of $\partial _{+\delta, \varepsilon}$ are, up to an  $\varepsilon$-error, equal the $B_x(\varepsilon)$-averages of these of     $\partial _{+\delta}$; hence, the mean curvatures of  $\partial V'=\partial _{+\delta, \varepsilon}$   are $>\phi(x)$ for sufficiently small $\varepsilon>0$.

Finally, in order to have $\partial V'$ inside rather than outside $V$, we apply the above to an interior equidistant hypersurface $\partial V_{-\delta}$ instead of  $\partial V=\partial V_0,$ where a minor readjustment of this argument is  needed  non-compact $\partial V$. 

\vspace {1mm}

This, together with   "essential example" allows  an alternative proof of

 \vspace {1mm}
 
  {\it   Smoothing of   Quasiregular Hypersurfaces }. Let $U$ be  an open domain
 in $X$ with quasiregular boundary and let the mean curvatures at all regular points of $\partial U$ are strictly minorized by a continous function $\psi$ on $X$, i.e. 
 $$mn.curv_x  (\partial U)>  \psi(x)  \mbox  { for all regular ponts } x\in \partial U.$$

{\it Then $U$ can be exhausted by closed subsets $U_i\subset U$ with smooth   boundaries $Y_i=\partial U_i$, where
the mean curvatures of these are  strictly  minorized by $\psi(x)$ at all $x\in \partial U_i$ and all $U_i$.}

 \vspace {1mm}

The two  basic examples where this smooth  approximation  is used in the present paper  are \vspace {1mm}

(1) {\it  Strictly Mean  Convex Bubbles $U\subset X$ with compact  boundaries $Y$.}

(2) {\it Minimal Hypersurfces $H$.}  \vspace {1mm}

In both cases   the  Almgren-Allard quasiregularity theorem applies and,  in the case (1), 
 allows a {\it smooth strictly mean convex approximation of $Y$  }  while in the case (2) one approximates the boundary $Y_\varepsilon$ of the 
 $\varepsilon$-neighbourhood $U\varepsilon(H)\subset X$ of $H$ by a {\it smooth 
 $c_\varepsilon$-concave hypersurface}, i.e. with $mn.curv (Y_\varepsilon) \leq c_\varepsilon$
 where $c_\varepsilon\to 0$ for $\varepsilon \to 0$.  \vspace {1mm}

 {\it Remarks.}(a)  When we    discussed  smoothing  minimal hypersurfaces $H$ with Joachim Lohkamp
  a few  years ago he, on one hand, 
 said he was well aware of possibility of such smoothing, but, on the other hand, he expressed a concern about singularities at  the focal points.

Focal points are invisible in the argument with bending and standard linear  smoothing (see section 3.4)  but the above  makes it clear why singularities at these points cause no additional complication. 
 
 (b)  The  corner  smoothing lemma. remain valid for the  $(n-k)$-mean 
convexity  for all $k$ but it fails, in general, for $k$-convexity, probably for all
 $k\neq 1, n-1$.
 
To see this  for even  $n-1\geq 4$ and $k=(n-1)/2\geq 2$, let $V\subset \mathbb R^{n-1}\subset \mathbb R^n$ be a compact domain with smooth boundary. Then there obviously exist  $C^\infty$-small perturbations  $V'$,   $V_+$
and $V_-$ of $V$ in $\mathbb R^n$, such that  $V_+$
and $V_-$ transversally meet along the boundary $\partial V'$ and bound together a 
domain $U'\subset \mathbb R^n$ that is $k$-convex away from the corner
along $\partial V'$. 

This $U'\supset V'$ can be seen as  a  small thickening of $V'$ that  is homeomorphic 
to $V'\times [0,1]=V\times [0,1]$.
Therefore, if the homology group  $H_{n-2}(V)\neq 0$, then    $H_{n-2}(U)\neq 0$ as well; hence,
 $U$ can not be approximated by smooth $k$-convex domains  if $k<n-2$.\vspace {2mm}

 \textbf {Acknowledgments.} I want to thank the anonymous referee who has indicated several errors in the original manuscript and has made useful suggestions.

  \end{document}